\documentclass[12pt]{amsart}
\usepackage{latexsym, amssymb, amsmath, amscd}
 \usepackage{eucal}

\theoremstyle{plain}
    \newtheorem{rema}{Remark}[section]
    \newtheorem{propo}[rema]{Proposition}

   \newtheorem{theo}[rema]{Theorem}
 \newtheorem{conj}[rema]{Conjecture}
   \newtheorem{defi}[rema]{Definition}
    \newtheorem{lemma}[rema]{Lemma}
    \newtheorem{corol}[rema]{Corollary}
     \newtheorem{exam}[rema]{Example}
  \newtheorem{rmk}[rema]{Remark}

	\newcommand{\nno}{\nonumber}

	\newcommand{\p}{\partial}

 \newcommand{\pf}{{\it Proof:}\hspace{2ex}}
 
 \newcommand{\epfv}{\hspace{1em}$\Box$\vspace{1em}}

\newcommand{\bR}{{\mathbb R}}
\newcommand{\bQ}{{\mathbb Q}}
\newcommand{\bN}{{\mathbb N}}

\newcommand{\cB}{{\mathcal B}}
\newcommand{\bC}{\mathbb C}

\newcommand{\cE}{{\mathcal E}}
\newcommand{\cR}{{\mathcal R}}
\newcommand{\cL}{{\mathcal L}}
\newcommand{\cT}{{\mathcal T}}
\newcommand{\cM}{{\mathcal M}}

\newcommand{\cA}{{\mathcal A}}

\newcommand{\cD}{{\mathcal D}}

\newcommand{\cx}{ \bC[\xi] }
\newcommand{\cz}{ \bC[z] }

\newcommand{\xiz}{{\xi, z}}
\newcommand{\cxz}{\bC[\xi, z]}

\newcommand{\im}{\rm{Im\,}}
\newcommand{\lin}{1\le i\le n}
\newcommand{\dlt}{\delta}

\renewcommand{\theequation}{\thesection.\arabic{equation}}
\renewcommand{\therema}{\thesection.\arabic{rema}}
\newcommand{\EAn}{\end{align*}}

\title[A Deformation of Commutative Polynomial Algebras]
{A Deformation of Commutative Polynomial Algebras in Even Numbers of Variables}

  \author{Wenhua Zhao}      

\begin{document}

\begin{abstract}
We introduce and study a deformation of commutative  polynomial algebras in even numbers of variables. We also discuss some connections and 
applications of this deformation to the generalized Laguerre orthogonal polynomials and 
the interchanges 
of right and left total symbols 
of differential operators of polynomial algebras.   
Furthermore, a more conceptual re-formulation 
for the {\it image} conjecture 
\cite{IC} is also given in terms of 
the deformed algebras. Consequently, 
the well-known {\it Jacobian} 
conjecture \cite{Ke} 
is reduced to an open  
problem on this deformation 
of polynomial algebras.
\end{abstract}

   

 \bibliographystyle{alpha}
    \maketitle


\renewcommand{\theequation}{\thesection.\arabic{equation}}
\renewcommand{\therema}{\thesection.\arabic{rema}}
\setcounter{equation}{0}
\setcounter{rema}{0}
\setcounter{section}{0}

\section{\bf  Introduction}\label{S1}

Let $\xi=(\xi_1, \xi_2, ..., \xi_n)$
and $z=(z_1, z_2, ..., z_n)$  
be $2n$ commutative free variables. 
Throughout this paper, we denote by 
$\cxz$, $\cz$ and $\cx$ the vector spaces 
(without any algebra structures)
over $\bC$  
of polynomials in $(\xiz)$, in $z$ and in $\xi$, respectively. 
The corresponding polynomial algebras will 
be denoted respectively by $\cA[\xi, z]$, 
$\cA[z]$ and $\cA[\xi]$. 

For any $\lin$, we set $\p_i\!:=\p_{z_i}$ and 
$\dlt_i=\p_{\xi_i}$. Denote by 
$\p=(\p_1, \p_2, ..., \p_n)$ 
and 
$\dlt=(\dlt_1, \dlt_2, ..., \dlt_n)$.
We also occasionally  use $\p_z$ and $\p_\xi$ 
to denote $\p$ and $\dlt$, respectively.

Set $\Omega\!:=\sum_{i=1}^n (\p_i \otimes \dlt_i +\dlt_i\otimes \p_i)$. For any $t\in \bC$, we define a new product $\ast_t$ for the vector space $\cxz$ by setting, 
for any $f, g\in \cxz$, 
\begin{align}\label{Def-ast-t}
f \ast_t g\!:= \mu\left(
e^{-t\, \Omega}(f\otimes g)\right),
\end{align} 
where $\mu:\cxz\otimes \cxz \to \cxz$ denotes 
the product map of the polynomial algebra 
$\cA[\xiz]$.

Denote by $\cB_t[\xiz]$ $(t\in \bC)$ 
the new algebra $(\cxz, \ast_t)$. 
For the case $t=1$, we also introduce 
the following short notation:
\begin{align}
\ast\!:&=\ast_{t=1}. \label{ast-t=1} \\
\cB[\xiz]\!:&=\cB_{t=1}[\xiz].\label{B-t=1}
\end{align}

Note that, when $t=0$, the algebra 
$\cB_{t=0}[\xiz]$ coincides with the 
usual polynomial algebra $\cA[\xiz]$.

In this paper, we first show that 
$\cB_t[\xiz]$ $(t\in \bC)$ gives 
a deformation of the polynomial 
algebra $\cA[\xiz]$. 
Actually, it is a trivial deformation 
in the sense of deformation theory. 
To be more precise, set 
\begin{align}
\Lambda\!:&=\sum_{i=1}^n \dlt_i\p_i. \label{Def-Lambda}\\
\Phi_t\!:&=e^{t\Lambda}=\sum_{m\ge 0} 
\frac {t^m \Lambda^m}{m!}. \label{Def-Phi-t}\\
\Phi\!:&=\Phi_{t=1}.\label{Def-Phi}
\end{align}
 
Note that, $\Phi_t$ for any $t\in \bC$ is  
a well-defined bijective linear map from 
$\cxz$ to $\cxz$, whose inverse map is 
given by $\Phi_{-t}=e^{-t\Lambda}$. 
This is because the differential 
operator $\Lambda$ of $\cxz$ is 
locally nilpotent, 
i.e.\@ for any $f(\xi, z)\in \cxz$, 
$\Lambda^m f(\xi, z)=0$ 
when $m\gg 0$.

With the notation fixed above, 
we will show that, for any $t\in \bC$,  
$\Phi_t:\cB_t[\xiz]\to \cA[\xiz]$ 
actually is an isomorphism of 
$\bC$-algebras (See Proposition \ref{GD-trivial} 
and Corollary \ref{D-trivial}).

Note that, from the point view of deformation theory, 
the deformation $\cB_t[\xiz]$ $(t\in \bC)$ 
is not interesting at all. But, surprisingly, 
as we will show in this paper, the algebra 
$\cB_t[\xiz]$ and the isomorphism $\Phi_t$ 
are actually closely related with the generalized Laguerre polynomials (See \cite{Sz}, \cite{PS} and \cite{AAR})
and the interchanges 
of right and left total symbols of differential 
operators of polynomial algebras. 

Furthermore, as we will show in Section \ref{S4}, 
the algebras $\cB_t[\xiz]$ $(t\in \bC)$ 
and the isomorphism $\Phi_t$ via their connections 
with the {\it image} conjecture proposed in 
\cite{IC} are also related with the {\it Jacobian} 
conjecture which was first proposed 
by O. H. Keller \cite{Ke} in $1939$ 
(See also \cite{BCW} and \cite{E}). 
Actually, the {\it Jacobian} conjecture can be 
viewed as a conjecture which, in some sense, 
just claims that the algebra $\cB_t[\xiz]$ $(t\ne 0)$  
should not differ or change too much from 
the polynomial algebra $\cA[\xiz]=\cB_{t=0}[\xiz]$. Therefore, from this point of view, the triviality of 
the deformation $\cB_t[\xiz]$ $(t\in \bC)$ 
(in the sense of deformation theory) can be 
viewed as a fact in favor of the 
{\it Jacobian} conjecture. 
For another interesting application 
of the isomorphism $\Phi$ to the 
{\it Jacobian} conjecture, see \cite{AGVC}. 

Considering the length of the paper, 
below we give a more detailed description 
for the contents and the arrangement 
of the paper.

In Subsection \ref{S2.1}, we prove some simple 
properties of the deformation $\cB_t[\xiz]$ $(t\in \bC)$ 
and the isomorphism $\Phi_t: \cB_t[\xiz]\to \cA[\xiz]$,  
which will be needed for the rest of this paper.  
In particular, in this subsection the triviality of the deformation $\cB_t[\xiz]$ $(t\in \bC)$ in the sense of deformation theory is proved in 
Proposition \ref{GD-trivial}. 
and Corollary \ref{D-trivial}.

In Subsection \ref{S2.2}, we show 
that, for different $t\in \bC$, the 
$\ell$-adic topologies 
induced by $\cB_t[\xiz]$ on the common base vector 
space $\cxz$ are different. But they are all 
homeomorphic to the $\ell$-adic topology  
induced by the polynomial algebra 
$\cA[\xiz]$ under the isomorphism 
$\Phi_t: \cB_t[\xiz] \to \cA[\xiz]$ 
(viewed as an automorphism of $\cxz$).
See Proposition \ref{L-topology} and also 
Corollary \ref{C2.10} for 
the precise statements.

In Subsection \ref{S2.3}, we study the induced 
isomorphism $(\Phi_t)_*$ $(t\in \bC)$ of 
$\Phi_t$ from the differential operator algebra, 
or the Weyl algebra $\cD_t[\xiz]$ of $\cB_t[\xiz]$
to the Weyl algebra $\cD[\xiz]$ of $\cA[\xiz]$. The main results 
of this subsection are 
Propositions \ref{Phi_*-P1} and \ref{Phi_*-P2}.
Proposition \ref{Phi_*-P1} says that 
the derivations $\p_{z_i}$ and 
$\p_{\xi_i}$ $(\lin)$ of $\cA[\xiz]$ are also 
derivations of $\cB_t[\xiz]$ for all $t\in \bC$ 
and are fixed by the isomorphism 
$(\Phi_t)_*$. Proposition \ref{Phi_*-P2} gives 
explicitly the images under $(\Phi_t)_*$ 
of the multiplication operators 
with respect to the product $\ast_t$ of 
$\cB_t[\xiz]$. 

In Section \ref{S3}, by using some results 
derived in Section \ref{S2}, we show in 
Theorem \ref{Interchange} 
that $\Phi=\Phi_{t=1}$ (resp.\@ $\Phi_{t=-1}$) 
as an automorphism of $\cxz$ actually 
coincides with the linear map 
which changes left (resp.\@ right) 
total symbols of differential operators 
of $\cA[z]$ to their right 
(resp.\@ left) total symbols.
Consequently, the products  
$\ast_{t=\pm 1}$ appear  
naturally when one derives  
left or right total symbols 
of certain differential operators 
of $\cA[z]$ (See Corollary \ref{C3.2}).
The results derived in this subsection 
also play some important roles in \cite{AGVC} 
in which among some other results 
a more straightforward proof for 
the equivalence of the {\it Jacobian} 
conjecture and the {\it vanishing} 
conjecture (See \cite{HNP} and \cite{GVC}) 
will be given. 
 
In Subsection \ref{S4.1}, we study the 
Taylor series expansion of 
polynomials in $\bC[\xiz]$ 
with respect to the new product 
$\ast_t$ and use it to give 
a more conceptual proof for 
the expansion of polynomials given 
in Eq.\,(\ref{D_t-Taylor-e2}).
This expansion was first proved 
in \cite{IC} and played a crucial role 
there in the proof of the implication of 
the {\it Jacobian} conjecture from the {\it image} conjecture (See Conjecture \ref{IC1}).

In Subsection \ref{S4.2}, we first recall 
the notion of the so-called {\it Mathieu} 
subspaces of commutative algebras
(See Definition \ref{M-ideal}), 
which was first introduced in \cite{GIC}, 
and also the {\it image} conjecture 
(See Conjecture \ref{IC1}) for 
the differential operators $\xi_i-t\p_i$ 
$(\lin)$ in terms of 
the notion of Mathieu subspaces.  
We then give a re-formulation of 
Conjecture \ref{IC1}  
in terms of the algebra 
$\cB_t[\xiz]$ $(t\in \bC)$ 
(See Conjecture \ref{IC2}) and 
show in Theorem \ref{IC1=IC2} that these two 
conjectures are equivalent to each other.
Since it has been shown in \cite{IC} that 
Conjecture \ref{IC1} implies {\it Jacobian} 
conjecture, hence so does 
Conjecture \ref{IC2}.

Consequently, via its connections with 
Conjecture \ref{IC2}, the {\it Jacobian} conjecture 
is reduced to an open problem on the deformation 
$\cB_t[\xiz]$ $(t\in \bC)$ of the polynomial algebra 
$\cA[\xiz]$. The open problem asks 
if the ideal $\xi\bC[\xiz]$ 
of $\cA[\xiz]$ generated by $\xi$ will 
remain to be a Mathieu subspace in 
the algebra $\cB_t[\xiz]$ for any $t\ne 0$. 
Note that any ideal is automatically 
a Mathieu subspace, but not conversely.
Therefore, the triviality 
(in the sense of deformation theory) 
of the deformation $\cB_t[\xiz]$ 
$(t\in \bC)$ proved in Proposition \ref{GD-trivial} 
can be viewed as a fact in favor of 
the {\it Jacobian} conjecture.

Section \ref{S5} is mainly on a connection  
of the algebra $\cB[\xiz]$, especially, 
its product $\ast$ with the multi-variable 
generalized Laguerre polynomials, 
and also some of the applications 
of this connection 
to both $\cB[\xiz]$ and the 
generalized Laguerre 
polynomials.

In Subsection \ref{S5.1}, we very briefly 
recall the definition of the 
(generalized) Laguerre polynomials 
$L_\alpha^{[\bf k]}(z)$ $({\bf k}, \alpha\in \bN^n)$  
(See Eqs.\,(\ref{Explicit-Formula})--(\ref{Mult-LPs}))  
and also the orthogonal property 
(See Theorem \ref{ortho}) of 
these polynomials.

In Subsection \ref{S5.2}, we show in 
Theorem \ref{ast2LP} that, for any 
${\bf k}, \alpha\in \bN^n$, we have 
\begin{align}\label{ast2LP-Int-e}
L_{\alpha}^{[{\bf k}]}( \xi z)
=\frac {(-1)^{|\alpha|}}{\alpha!} \, 
 \xi^{-\bf k}(\xi^{\alpha+{\bf k}} \ast z^\alpha) 
=\frac {(-1)^{|\alpha|}}{\alpha!} \, 
z^{-{\bf k}}(\xi^\alpha \ast z^{\alpha+{\bf k}}).
\end{align}

Consequently, the generalized Laguerre 
polynomials $L_\alpha^{[\bf k]}(z)$ 
$({\bf k}, \alpha\in \bN^n)$ 
can be obtained by evaluating the polynomials 
$ \xi^{-\bf k}(\xi^{\alpha+{\bf k}} \ast z^\alpha)$ or 
$z^{-{\bf k}}(\xi^\alpha \ast z^{\alpha+{\bf k}})$ at 
$\xi=(1, 1,...,1)$. Note that the evaluation map at 
$\xi=(1, 1,...,1)$ is not an algebra homomorphism 
from $\cB[\xiz]$ to $\cA[\xiz]$. Otherwise,  
the generalized Laguerre polynomials 
would be trivialized.

In the first part of Subsection \ref{S5.3}, 
we use certain results of the generalized Laguerre 
polynomials and the connection 
in Eq.\,(\ref{ast2LP-Int-e}) above 
to derive more properties 
on the polynomials $\xi^\alpha \ast z^\alpha$ which, 
by Proposition \ref{freeness}, $(c)$, are actually 
the monomials of $\xi$ and $z$ in the 
new algebra $\cB[\xiz]$. 

For example, by using the connection  
in Eq.\,(\ref{ast2LP-Int-e}) and 
I. Schur's irreducibility theorem \cite{S1} of 
the Laguerre polynomials in one variable, 
we immediately have that, when $n=1$, 
the monomials $\xi^m\ast z^m$ 
$(m\ge 2)$ of $\cB[\xiz]$ are 
actually irreducible over $\bQ$ 
(See Theorem \ref{Irreduce}). 
Furthermore, by using I. Schur's irreducibility theorem \cite{S2} and 
M. Filaseta and T.-Y. Lam's irreducibility 
theorem \cite{FL} on the generalized  
Laguerre polynomials, 
we have that, all but finitely many of 
the polynomials $\xi^{-k}(\xi^{m+k}\ast z^m)$ and   
$z^{-k}(\xi^m\ast z^{m+k})$ $(m, k\in \bN)$
are irreducible over $\bQ$ 
(See Theorem \ref{Irreduce-2}).

In the second part of Subsection \ref{S5.3},
we use the connection given in Eq.\,(\ref{ast2LP-Int-e}) 
and certain results of $\cB[\xiz]$ derived 
in Section \ref{S2} to give new proofs, first, 
for some recurrent formulas of the 
generalized Laguerre polynomials 
(See Proposition \ref{Recur-Formula}) 
and, second, for the fact that the 
generalized Laguerre polynomials satisfy  
the so-called 
{\it associate Laguerre differential equation} 
(See Theorem \ref{Aso-Lag-Equation}). 
At the end of this subsection,  
we draw the reader's attention to a conjecture, 
Conjecture \ref{Conj-Lgr}, on the generalized 
Laguerre polynomials, which is still open 
even for the classical Laguerre polynomials 
in one variable.

{\bf Acknowledgment} The author would like to thank 
the anonymous referees for pointing out many  
typos, minor errors of the previous version of this paper, 
and also for suggesting the new  
proof of Lemma \ref{Irreduce-lemma} without the condition 
that the base filed $K$ has infinitely many elements.

\renewcommand{\theequation}{\thesection.\arabic{equation}}
\renewcommand{\therema}{\thesection.\arabic{rema}}
\setcounter{equation}{0}
\setcounter{rema}{0}

\section{\bf The Deformation $\cB_t[\xiz]$ of 
the Polynomial Algebra $\cA[\xiz]$} \label{S2}

In this section, we first derive in Subsection \ref{S2.1} 
some properties and identities for the algebra $\cB_t[\xiz]$ 
$(t\in \bC)$. In Subsection \ref{S2.2}, we show that, 
for different $t\in \bC$, the $\ell$-adic topologies 
induced by the algebras  $\cB_t[\xiz]$ $(t\in \bC)$ 
on the common base vector space $\cxz$ 
are different. But they are all homeomorphic 
under the isomorphism 
$\Phi_t:\cB_t[\xiz]\to \cA[\xiz]$
to the $\ell$-adic topology on $\cxz$ induced 
by the polynomial algebra $\cA[\xiz]$ 
(See Proposition \ref{L-topology} and also 
Corollary \ref{C2.10}).

In Subsection \ref{S2.3}, we study 
the isomorphism $(\Phi_t)_*$ induced by 
$\Phi_t$ from the Weyl algebra of $\cB_t[\xiz]$ 
to the Weyl algebra of $\cA[\xiz]$. 
The main results in this subsection are 
Propositions \ref{Phi_*-P1} and \ref{Phi_*-P2}.

\subsection{Some Properties of the Algebras $\cB_t[\xiz]$} 
\label{S2.1}
First, one remark on notation and convention is that, we will freely use throughout this paper some commonly used 
multi-index notations and conventions.
For instance, for $n$-tuples  
$\alpha=(k_1, k_2, ..., k_n)$ and $\beta=(m_1, m_2, ..., m_n)$ of non-negative integers, we have
\allowdisplaybreaks{
\begin{align*}
|\alpha|&=\sum_{i=1}^n k_i. \\
\alpha! &=k_1!k_2!\cdots k_n!. \\
\binom{\alpha}{\beta} &=
\begin{cases} \frac{\alpha!}{\beta! (\alpha-\beta)!} 
&\mbox{if  $k_i\ge m_i$ for all $1\le i\le n$}; \\
0, &\mbox{otherwise.}
\end{cases}
\end{align*} }

The notation and convention fixed in the previous 
section will also be used throughout this paper. 
 
The first main result of this section 
is the following proposition.

\begin{propo}\label{GD-trivial}
For any $t\in \bC$ and 
$f, g\in \cxz$, we have
\begin{align}\label{GD-trivial-e1}
\Phi_t (f \ast_t g) = \Phi_t (f) \Phi_t (g).
\end{align}
\end{propo}

\pf We first set 
\begin{align}\label{GD-trivial-pe1}
f \ast_t' g \!:= \Phi_t^{-1} \big(\Phi_t (f) \Phi_t (g)\big)=\Phi_{-t} \big(\Phi_t (f) \Phi_t (g)\big).
\end{align}
for any $f, g\in \bC[\xiz]$.

We view $t$ as a formal parameter which commutes with $\xi$ and $z$. 
Then, by Eqs.\,(\ref{Def-ast-t}), (\ref{GD-trivial-pe1}) and 
the fact that the differential operators $\Lambda$ and $\Omega$ are locally nilpotent on $\cxz$ and $\cxz\otimes\cxz$, respectively, we see that $f \ast_t g$ and   
$f\ast_t' g$ are polynomials in $t$ 
with coefficients in $\bC[\xiz]$. 
Furthermore, by setting $t=0$ in 
Eqs.\,(\ref{Def-ast-t}) and (\ref{GD-trivial-pe1}), 
we see that the constant terms (with respect to $t$) 
of $f \ast_t g$ and $f\ast_t' g$ are both 
$fg\in \bC[\xiz]$. In other words, we have
\begin{align}\label{GD-trivial-pe2}
f \ast_t g\left. \right|_{t=0}=f \ast_t' g
\left. \right|_{t=0}=fg.
\end{align}

From Eq.\,(\ref{Def-ast-t}),  we have, 
\begin{align}
\frac{\p}{\p t}(f\ast_t g)&= 
-\mu\left(\, e^{-t\Omega} (\Omega(f\otimes g)) 
\, \right)\label{GD-trivial-pe3} \\
&= -\sum_{i=1}^n \mu \left( e^{-t\Omega} \left
( (\dlt_i f) \otimes (\p_i g) + (\p_i f) \otimes 
(\dlt_i g)  \right) \right) \nno \\
&= -\sum_{i=1}^n \left( (\dlt_i f) \ast_t 
(\p_i g) + (\p_i f) \ast_t (\dlt_i g) \right ). \nno 
\end{align}
 
On the other hand, from Eq.\,(\ref{GD-trivial-pe1}), we have, 
\begin{align*}
\frac{\p}{\p t}(f&\ast_t' g) = 
\frac{\p}{\p t} \left( e^{-t\Lambda}
( (e^{t\Lambda}f)\,(e^{t\Lambda}g) ) \right)\\
&= e^{-t\Lambda} \left( 
-\Lambda ( (e^{t\Lambda}f)\,(e^{t\Lambda}g) ) 
+ (e^{t\Lambda}\Lambda f)\,(e^{t\Lambda}g) +
(e^{t\Lambda}f)\,(e^{t\Lambda}\Lambda g)  \right).
\end{align*}

Note that, for any $u, v\in \cxz$, 
it is easy to check that we have the following identity:
\begin{align*}
\Lambda(uv)
=(\Lambda u)v+u(\Lambda v)+\sum_{i=1}^n 
\left( \, (\dlt_iu) (\p_i v)
+(\p_i u) (\dlt_i v) \, \right).
\end{align*}

By the last two equations above and 
also Eq.\,(\ref{GD-trivial-pe1}), we have
\begin{align}
\frac{\p}{\p t}(f\ast_t' g)
&=-\sum_{i=1}^n e^{-t\Lambda} \left( \, 
( (e^{t\Lambda}\dlt_i f)\,(e^{t\Lambda}\p_i g) ) 
+ (e^{t\Lambda}\p_i f)\,(e^{t\Lambda}\dlt_ig)\,  \right) 
\label{GD-trivial-pe4}\\
&=-\sum_{i=1}^n\left( \, (\dlt_if) \ast_t' 
(\p_i g) + (\p_i f) \ast_t' (\dlt_i g) \, \right). \nno 
\end{align}

Next, we use the induction on 
$(\deg f+\deg g)$ to show Eq.\,(\ref{GD-trivial-e1}).
First, when $\deg f+\deg g=0$, i.e.\@  
both $f$ and $g$ have degree zero,  
it is easy to see from Eqs.\,(\ref{Def-ast-t}) 
and (\ref{GD-trivial-pe1}) that  
$f\ast_t g = f\ast_t' g=fg$ 
in this case. 

In general, by 
Eqs.\,(\ref{GD-trivial-pe3}), (\ref{GD-trivial-pe4}) 
and also the induction assumption, we have 
\begin{align}\label{GD-trivial-pe5}
\frac{\p}{\p t}(f\ast_t g)=\frac{\p}{\p t}(f \ast_t' g).
\end{align}

Since $f\ast_t g$ and $f\ast_t' g$ are polynomials in $t$ 
with coefficients in $\bC[\xiz]$ and both satisfy  
Eqs.\,(\ref{GD-trivial-pe2}) and (\ref{GD-trivial-pe5}), 
it is easy to see that they must 
be equal to each other. Hence, 
Eq.\,(\ref{GD-trivial-e1}) holds.
\epfv

\begin{corol}\label{D-trivial}
For any $t\in \bC$, $\Phi_t: \cB_t[\xiz] \to \cA[\xiz]$ 
is an isomorphism of algebras. 
Therefore, in the sense of deformation theory, 
the deformation $\cB_t[\xiz]$ is a trivial 
deformation of the commutative polynomial algebra $\cA[\xiz]$. 
\end{corol}

Next we derive some properties of the algebras 
$\cB_t[\xiz]$ $(t\in \bC)$, which will be needed 
for the rest of this paper.

\begin{lemma}\label{B-Lemma}
For any $f, g\in \bC[\xiz]$, we have
\begin{align}\label{B-Lemma-e1}
f\ast_t g=\sum_{\alpha, \beta\in \bN^n} 
\frac{(-t)^{|\alpha|+|\beta|}}{\alpha!\beta!} \,\,
(\dlt^\beta \p^\alpha f) (\p^\beta \dlt^\alpha g).
\end{align}
\end{lemma}

\pf Note first that, for any $1\le i, j\le n$, 
$\p_i\otimes \dlt_i$ and $\dlt_j\otimes \p_j$ 
commute with each other. So we have 
\begin{align}
e^{-t\, \Omega}&=e^{-t\sum_{i=1}^n \dlt_i\otimes \p_i}\, 
e^{-t\sum_{i=1}^n \p_i\otimes \dlt_i}, \\
e^{-t\sum_{i=1}^n \p_i\otimes \dlt_i}
&=\prod_{i=1}^n e^{-t (\p_i\otimes \dlt_i)}
=\prod_{i=1}^n\sum_{k_i\ge 0} 
\frac{(-t)^{k_i}}{k_i!}\, 
(\p_i^{k_i}\otimes \dlt_i^{k_i}) \\
&=\sum_{\alpha\in \bN^n}
\frac{(-t)^{|\alpha|}}{\alpha!} \,\,
(\p^\alpha \otimes \dlt^\alpha). \nno
\end{align}

Similarly, 
\begin{align}
e^{-t\sum_{i=1}^n \dlt_i \otimes \p_i}
=\sum_{\beta \in \bN^n} 
\frac{(-t)^{|\beta|}}{\beta!} \,\,
(\dlt^\beta \otimes \p^\beta).
\end{align}

Then it is easy to see that Eq.\,(\ref{B-Lemma-e1}) follows directly 
from Eq.\,(\ref{Def-ast-t}) and the last three equations above.
\epfv

\begin{propo}\label{B-Propo}
$(a)$ For any $\lambda_i(\xi)\in \bC[\xi]$, 
$p_i(z)\in \bC[z]$ $(i=1, 2)$, we have 
\begin{align}
\lambda_1 (\xi) \ast_t \lambda_2 (\xi)&= 
\lambda_1 (\xi)  \lambda_2 (\xi). \label{B-Propo-e1} \\
p_1(z)\ast_t p_2(z)&=p_1(z) p_2(z).\label{B-Propo-e2} \\
\lambda(\xi) \ast_t p(z)=
\sum_{\alpha\in \bN^n} &
\frac {(-1)^{|\alpha|}t^{|\alpha|}}{\alpha!} 
(\dlt^\alpha \lambda(\xi))(\p^\alpha p(z)).
\label{B-Propo-e3} 
\end{align}

$(b)$ For any $\lambda(\xi)\in \bC[\xi]$, $p(z)\in \bC[z]$ and 
$g(\xiz)\in \bC[\xiz]$, we have
\begin{align}
\lambda(\xi) \ast_t g(\xiz)&= \lambda(\xi-t\p)g(\xiz). 
\label{B-Propo-e4}\\
p(z)\ast_t g(\xiz) &=p(z-t\dlt)g(\xiz). \label{B-Propo-e5}
\end{align}
\end{propo}

Note that the components $\xi_i-t\p_i$ $(\lin)$ of 
the $n$-tuple $\xi-t\p$ in Eq.\,(\ref{B-Propo-e4}) 
commute with one another. So the substitution 
$\lambda(\xi-t\p)$ of $\xi-t\p$ into the polynomial 
$\lambda(\xi)$ is well-defined. Similarly, 
the substitution $p(z-t\dlt)$ in 
Eq.\,(\ref{B-Propo-e5}) is 
also well-defined.\\

\pf Eqs.\,(\ref{B-Propo-e1})--(\ref{B-Propo-e3}) 
follow directly from 
Eq.\,(\ref{B-Lemma-e1}).

To show Eq.\,(\ref{B-Propo-e4}), first,  
by Eq.\,(\ref{B-Lemma-e1}), we have
\begin{align}
\lambda(\xi) \ast_t g(\xiz)=
\sum_{\alpha\in \bN^n} &
\frac {(-1)^{|\alpha|}t^{|\alpha|}}{\alpha!} 
(\dlt^\alpha \lambda(\xi)) (\p^\alpha g(\xiz)).
\label{B-Propo-pe1} 
\end{align}

Second, note that the multiplication operators 
by $\xi_i$ $(\lin)$ and the derivations  
$\p_j$ $(1\le j\le n)$ commute. By using the Taylor 
series expansion of $\lambda(\xi-t\p)$ at $\xi$, 
we have 
\begin{align}
\lambda(\xi-t\p )g(\xiz)&=
\left( \sum_{\alpha\in \bN^n} 
\frac {(-1)^{|\alpha|}t^{|\alpha|}}{\alpha!} 
(\dlt^\alpha \lambda)(\xi) \p^\alpha \right) g(\xiz) 
\label{B-Propo-pe2}\\
&=\sum_{\alpha\in \bN^n} 
\frac {(-1)^{|\alpha|}t^{|\alpha|}}{\alpha!} 
(\dlt^\alpha \lambda(\xi)) (\p^\alpha g(\xiz)).\nno
\end{align}

Hence, Eq.\,(\ref{B-Propo-e4}) follows 
from the last two equations.  Eq.\,(\ref{B-Propo-e5}) 
can be proved similarly.
\epfv

\begin{lemma}\label{L2.1}
For any $t\in \bC$, $\lambda(\xi)\in \bC[\xi]$ 
and $p(z)\in \bC[z]$, we have 
\allowdisplaybreaks{
\begin{align}
\Phi_t (\lambda(\xi))&=\lambda(\xi). \label{L2.1-e1}\\
\Phi_t (p(z))&=p(z). \label{L2.1-e2} \\
\Phi_t (\lambda(\xi) p(z)) 
&=\lambda(\xi) \ast_{-t} p(z).   \label{L2.1-e3} 
\end{align}}
\end{lemma}

\pf Since $\Lambda (\lambda(\xi))=\Lambda(p(z))=0$,  
$\Phi_t=e^{t\Lambda}$ fixes $\lambda(\xi)$ and $p(z)$. 
Hence we have Eqs.\,(\ref{L2.1-e1}) and (\ref{L2.1-e2}).

To show Eq.\,(\ref{L2.1-e3}), 
by Eqs.\,(\ref{GD-trivial-e1}), (\ref{L2.1-e1}) and (\ref{L2.1-e2}), 
we have   
\begin{align*}
\lambda(\xi) \ast_t p(z) &= 
\Phi_{-t}\left(\, \Phi_t(\lambda(\xi)) \Phi_t(p(z)) \, \right)  \\
&= \Phi_{-t}( \lambda(\xi)p(z) ).
\end{align*}
Replacing $t$ be $-t$ in the equation above, 
we get Eq.\,(\ref{L2.1-e3}).
\epfv

\begin{propo}\label{freeness}
For any $t\in \bC$, the following statements hold. 

$(a)$ The subspaces 
$\bC[\xi]$ and $\bC[z]$ of $\cB_t[\xiz]$
are closed under the product $\ast_t$ and 
hence, are actually subalgebras of 
$\cB_t[\xiz]$.

$(b)$ As associative algebras, $(\bC[\xi], \ast_t)$ 
and $(\bC[z], \ast_t)$ are identical as 
the usual polynomial algebras $\cA[\xi]$ 
and $\cA[z]$ in $\xi$ and $z$, respectively.

$(c)$ $\cB_t[\xiz]$ is a commutative free 
algebra generated freely by $\xi_i$ 
and $z_i$ $(\lin)$. The set of the monomials generated 
by $\xi_i$ and $z_i$ $(\lin)$ in $\cB_t[\xiz]$ 
is given by 
$\{\xi^\alpha\ast_t z^\beta\,|\, \alpha, \beta\in \bN^n\}$.
\end{propo}

\pf Note that $(a)$ and $(b)$ follow immediately from 
Eqs.\,(\ref{B-Propo-e1}) and (\ref{B-Propo-e2}).

To show $(c)$, first, by Eqs.\,(\ref{L2.1-e1}) 
and (\ref{L2.1-e2}), we know that 
the algebra isomorphism 
$\Phi_{-t}=\Phi_t^{-1}: \cA[\xiz]\to \cB_t[\xiz]$ 
as a linear map from $\cxz$  to $\cxz$ 
fixes $\xi_i$ and $z_i$ $(\lin)$. Hence, 
$\cB_t[\xiz]$ is a commutative free 
algebra generated freely by $\xi_i$ 
and $z_i$ $(\lin)$.

The second part of $(c)$ follows from 
Eqs.\,(\ref{B-Propo-e1}), (\ref{B-Propo-e2}) 
and the fact that the product $\ast_t$ 
is associative and commutative.
\epfv

The next two lemmas will be needed 
in Subsection \ref{S5.3}.

\begin{lemma}\label{monomial-eta}
For any $t\in \bC$ and $\alpha, \beta\in \bN$, 
\begin{align}\label{monomial-eta-e1}
(z\p-\xi\dlt)(\xi^\alpha \ast_t z^\beta)= 
(|\beta|-|\alpha|)(\xi^\alpha \ast_t z^\beta),
\end{align}
where $z\p-\xi\dlt\!:=\sum_{i=1}^n (z_i\p_i-\xi_i \dlt_i)$.
\end{lemma}

\pf First, by Euler's lemma, we have 
\begin{align}\label{monomial-eta-pe1}
(z\p-\xi\dlt)(\xi^\alpha z^\beta)= 
(|\beta|-|\alpha|)(\xi^\alpha z^\beta).
\end{align}

Second, note that $z\p-\xi\dlt$ commutes with $\Lambda$, 
hence also with $\Phi_t$ for any $t\in \bC$. 
Apply $\Phi_{-t}$ to Eq.\,(\ref{monomial-eta-pe1}), 
we get
\begin{align*}
(z\p-\xi\dlt)\Phi_{-t}(\xi^\alpha z^\beta)= 
(|\beta|-|\alpha|)\Phi_{-t}(\xi^\alpha z^\beta).
\end{align*}
Then, by Eq.\,(\ref{L2.1-e3}) with $t$ replaced by $-t$, 
Eq.\,(\ref{monomial-eta-e1}) follows 
from the equation above.
\epfv

\begin{lemma}\label{Symmetry}
For any $\lambda_i(\xi)\in \bC[\xi]$ 
and $p_i(z)\in \bC[z]$ $(i=1, 2)$, we have
\begin{align}\label{Symmetry-e1}
(\lambda_1(\xi) p_1(z)) \ast_t 
(\lambda_2(\xi) p_2(z))= 
(\lambda_1(\xi) \ast_{t} p_2(z) )\, 
(\lambda_2(\xi)\ast_{t} p_1(z)).
\end{align}
\end{lemma}

\pf First, by Eq.\,(\ref{B-Lemma-e1}), 
we have 
\allowdisplaybreaks{
\begin{align*}
&\quad (\lambda_1(\xi) p_1(z)) \ast_t 
(\lambda_2(\xi) p_2(z)) \\
&=\sum_{\alpha, \beta \in \bN^n} 
\frac {(-t)^{|\alpha|+|\beta|}}{\alpha!\beta!} 
\left( ( \dlt^\alpha \lambda_1(\xi)) (\p^\beta p_1(z)) \right)
\left( (\dlt^\beta \lambda_2(\xi)) (\p^\alpha p_2(z)) \right) \\
\intertext{Taking sum over $\alpha\in \bN^n$ 
and applying Eq.\,(\ref{B-Propo-e3}):}
&= (\lambda_1(\xi) \ast_{t} p_2(z) )
\sum_{\beta \in \bN^n} 
\frac {(-t)^{|\beta|}}{\beta!} 
 (\p^\beta p_1(z))(\dlt^\beta \lambda_2(\xi))  \\
\intertext{Taking sum over $\beta\in \bN^n$ 
and applying Eq.\,(\ref{B-Propo-e3}):}
&=(\lambda_1(\xi) \ast_{t} p_2(z))
(\lambda_2(\xi)\ast_{t} p_1(z)).
\end{align*} }
Hence we get Eq.\,(\ref{Symmetry-e1}).
\epfv

\subsection{The $\ell$-adic Topologies Induced by $\cB_t[\xiz]$ on $\cxz$}\label{S2.2}

We have seen that the algebras $\cB_t[\xiz]$ $(t\in \bC)$ 
share the same base vector space $\cxz$ and, 
by Proposition \ref{freeness}, $(c)$, 
they are all commutative free algebras generated 
freely by $\xi$ and $z$. Therefore, we may talk about 
the $\ell$-adic topologies on $\cxz$ induced by 
the algebras $\cB_t[\xiz]$ $(t\in \bC)$, 
which are defined as follows.

For any $t\in \bC[\xiz]$ and $m\ge 0$, 
set $U_{t, m}$ to be the subspace of $\cxz$ 
spanned by the monomials $\xi^\alpha \ast_t z^\beta$ 
of $\cB_t[\xiz]$ with $\alpha, \beta\in \bN^n$ 
and $|\alpha+\beta|\ge m$.
The $\ell$-adic topology 
induced from the algebra $\cB_t[\xiz]$ is the 
topology whose open subsets are the subsets 
generated by $U_{t, m}$ $(m\in \bN)$ and 
their translations by elements 
of $\cB_t[\xiz]$. We denote by $\cT_t$ 
this topology on $\cxz$.

The main result of this subsection is the following proposition.

\begin{propo}\label{L-topology}
$(a)$ For any $s\ne t\in \bC$, we have $\cT_s \ne \cT_t$.

$(b)$ For any $t\in \bC$, the algebra isomorphism 
$\Phi_t: (\cB_t[\xiz], \cT_t) \to (\cA[\xiz], \cT_0)$ 
is also a homeomorphism of topological spaces. 
Consequently, $(\cB_t[\xiz], \cT_t)$ $(t\in \bC)$ 
as topological spaces are all homeomorphic.
\end{propo}

\pf $(a)$ Let $\{\alpha_m \in \bN^n\,|\, m\ge 1\}$ be 
any sequence of elements of $\bN^n$ such that 
$|\alpha_m|=m$ for any $m\ge 1$.

Set $u_m\!:=\xi^{\alpha_m}\ast_t z^{\alpha_m}$ 
for any $m\ge 1$. Then, by the definition of 
$\cT_t$, we see that the sequence $\{u_m\}$ 
converges to $0 \in \cxz$ with respect 
to the topology $\cT_t$.
  
But, on the other hand, set 
$r:=s-t\ne 0$.  Then, by Eq.\,(\ref{B-Propo-e4}), we have
\begin{align*}
u_m&=\xi^{\alpha_m}\ast_t z^{\alpha_m}=
(\xi-t\p)^{\alpha_m}  z^{\alpha_m} 
=((\xi-s\p)+r\p)^{\alpha_m} z^{\alpha_m} \\
&= \sum_{\substack{\beta, \gamma \in \bN^n \\
\beta+\gamma=\alpha_m}} 
\binom {\alpha_m}{\gamma} (\xi-s\p)^\gamma (\p^\beta z^{\alpha_m}) \\
&= \sum_{\substack{\beta, \gamma \in \bN^n \\
\beta+\gamma=\alpha_m}} 
\binom {\alpha_m}{\gamma} \xi^\gamma \ast_s (\p^\beta z^{\alpha_m}) 
\equiv \alpha_m! \mod (U_{s, 0}).
\end{align*}

From the equation above, we see that 
the sequence $\{u_m\}$ does not converge 
to $0\in \cxz$ with respect 
to the topology $\cT_s$. Hence 
$\cT_s \ne \cT_t$.

$(b)$ Note that $\cB_{t=0}[\xiz]$ is the usual 
polynomial algebra $\cA[\xiz]$ and 
$\Phi_t: \cB_t[\xiz] \to \cA[\xiz]$ is an algebra 
isomorphism. Furthermore, from   
Eqs.\,(\ref{GD-trivial-e1}), (\ref{L2.1-e1}) and (\ref{L2.1-e2}), 
we have 
\begin{align}\label{Phi-M2M}
\Phi_t(\xi^\alpha\ast_t z^\beta)=\xi^\alpha z^\beta
\end{align}
for any $\alpha, \beta \in \bN^n$.

Therefore, for any $m\ge 0$, 
we have, $\Phi_t (U_{t, m})=U_{0, m}$ 
and $\Phi_t^{-1}(U_{0, m})=U_{t, m}$. 
Hence, we have $(b)$.
\epfv

Actually, the proof above also shows that 
Proposition \ref{L-topology} 
also holds for the following topologies on $\cxz$
induced by the free algebras $\cB_t[\xiz]$ $(t\in \bC)$.

For any $t\in \bC[\xiz]$ and $m\ge 0$, set 
\begin{align}\label{Um'} 
U_{t, m}'\!:=
\sum_{\substack{\alpha \in \bN^n;\, \\|\alpha|\ge m}} 
\xi^\alpha\ast_t\cz.
\end{align}

Denote by $\cT'_t$ the topology on $\cxz$ generated by 
$U_{t, m}'$ and their translations (as open subsets). 
Then, by a similar argument as in the proof of  
Proposition \ref{L-topology}, it is easy to see  
that the following corollary also holds.

\begin{corol}\label{C2.10} 
$(a)$ For any $s\ne t\in \bC$, we have $\cT_s{}' \ne \cT_t{}'$.

$(b)$ For any $t\in \bC$, the algebra isomorphism 
$\Phi_t: (\cB_t[\xiz], \cT_t{}') \to (\cA[\xiz], \cT_0{}')$ 
is also a homeomorphism of topological spaces. 
\end{corol}

Note that, due to the symmetric roles played by $\xi$ and $z$, 
the corollary above also holds if $\xi$ in 
Eq.\,(\ref{Um'}) is replaced by $z$.

\subsection{The Induced Isomorphism $(\Phi_t)_\ast$ 
on Differential Operator Algebras}\label{S2.3}

For any $t\in \bC$, denote by 
$\cD_t[\xiz]$ the differential 
operator algebra or the Weyl 
algebra of $\cB_t[\xiz]$, i.e.\@ the associative 
algebra generated by the $\bC$-derivations and the multiplication operators of 
the algebra $\cB_t[\xiz]$. 
Since 
$\Phi_t: \cB_t[\xiz]\to \cA[\xiz]$
is an algebra isomorphism 
(See Corollary \ref{D-trivial}), 
it induces an algebra isomorphism,   
denoted by 
$(\Phi_t)_*: \cD_t[\xiz] \to \cD[\xiz]$,
from the Weyl algebra  $\cD_t[\xiz]$ of $\cB_t[\xiz]$
to the Weyl algebra  $\cD[\xiz]$
of $\cA[\xiz]$. 

Recall that the induced map 
$(\Phi_t)_*$ is defined by setting 
\begin{align}\label{Def-Phi-*}
(\Phi_t)_*(\psi)=\Phi_t\circ \psi \circ \Phi_t^{-1}
=\Phi_t\circ \psi \circ \Phi_{-t} 
\end{align}
for any $\psi\in \cD_t[\xiz]$.
   
The main result of this subsection are 
the following two propositions, 
even though their proofs 
are very simple.

\begin{propo}\label{Phi_*-P1}
For any $t\in \bC$, the following statements hold.

$(a)$ $\p_i$ and $\dlt_i$ $(\lin)$ are also derivations 
of $\cB_t[\xiz]$.

$(b)$ For any $\lin$, we have 
\begin{align}
(\Phi_t)_*(\p_i)&=\p_i, \label{Phi_*-P1-e1}\\
(\Phi_t)_*(\dlt_i)&=\dlt_i,\label{Phi_*-P1-e2}
\end{align}
\end{propo}

\pf Note first that $\p_i$ and $\dlt_i$ $(\lin)$ commute 
with $\Lambda$, hence also with $\Phi_t$ 
for any $t\in \bC$. Then, Eqs.\,(\ref{Phi_*-P1-e1}) 
and (\ref{Phi_*-P1-e2}) follows immediately 
from this fact and the definition of $(\Phi_t)_*$ 
given in Eq.\,(\ref{Def-Phi-*}).

$(a)$ follows from the general fact that the induced map of any algebra isomorphism maps derivations 
to derivations. 
It can also be checked  
directly as follows.

For any $f, g\in \cB_t[\xiz]$, 
by Eq.\,(\ref{GD-trivial-e1}) and the fact that 
$\p_i$ $(\lin)$ commute with $\Phi_{t}$ $(t\in \bC)$,  
we have  
\begin{align*}
\p_i (f\ast_t g)&=\p_i \Big(\Phi_{-t} \big(\Phi_t(f)\Phi_t(g)\big) \Big)
=\Phi_{-t} \Big(\p_i\big( \Phi_t(f)\Phi_t(g)\big) \Big) \nno \\
&= \Phi_{-t} \big( (\p_i \Phi_t(f) ) \Phi_t(g)\big)  
+\Phi_{-t} \big( \Phi_t(f) (\p_i \Phi_t(g) ) \big)\\
&= \Phi_{-t} \big( \Phi_t(\p_i f)  \Phi_t(g)\big)  
+\Phi_{-t} \big( \Phi_t(f) \Phi_t(\p_i g)  \big)\\
&=(\p_i f)\ast_t g + f\ast_t(\p_i g). \nno 
\end{align*}

Similarly, we can show that $\dlt_i$ $(\lin)$ 
are also derivations of $\cB_t[\xiz]$.
\epfv

\begin{corol}\label{ast-Leibniz}
For any $\alpha, \beta, \gamma\in \bN^n$, we have
\begin{align}
\p^\gamma(\xi^\alpha \ast_t z^\beta)&=\gamma!\,
\binom\beta\gamma \, (\xi^\alpha \ast_t z^{\beta-\gamma}),
\label{ast-Leibniz-e1} \\
\dlt^\gamma(\xi^\alpha \ast_t z^\beta)&=\gamma!\,
\binom\alpha\gamma \, (\xi^{\alpha-\gamma}\ast_t z^\beta).
\label{ast-Leibniz-e2}
\end{align}
\end{corol}

\pf Note that, by Eqs.\,(\ref{B-Propo-e1}) 
and (\ref{B-Propo-e2}), we know that, 
for any $\alpha, \beta\in \bN^n$,  
$\xi^\alpha \ast_t z^\beta$ 
will remain the same if we replace the (usual) product 
of $\cA[\xiz]$ in the factors $\xi^\alpha$ and 
$z^\beta$ by the product $\ast_t$ of $\cB_t[\xiz]$. 
By Proposition \ref{Phi_*-P1}, $(a)$, we know that 
$\p_i$ and $\dlt_i$ $(\lin)$ are also 
the derivations of $\cB_t[\xiz]$. 
From these two facts, it is easy to see 
that both equations in the corollary  
hold.
\epfv

\begin{propo}\label{Phi_*-P2}
For any $t\in \bC$ and   
$f(\xi, z)\in \bC[\xiz]$, $(\Phi_t)_*$ maps 
the multiplication operator of $\cB_t[\xiz]$   
by $f(\xi, z)$ 
$($with respect to the product $\ast_t$$)$ 
to the multiplication operator of $\cA[\xiz]$
by $\Phi_t(f(\xi, z))$ 
$($with respect to 
the product of $\cA[\xiz]$$)$.   
\end{propo}

\pf We denote by $\psi_f$ 
the multiplication operator of $\cB_t[\xiz]$   
by $f(\xi, z)$ 
$($with respect to the product $\ast_t$$)$.  
Then for any $u(\xiz)\in \bC[\xiz]$, by  
Eqs.\,(\ref{Def-Phi-*}) and (\ref{GD-trivial-e1})  
we have 
\begin{align*} 
(\Phi_t)_*(\psi_f)u(\xiz)&=
(\Phi_t\circ \psi_f \circ \Phi_t^{-1})u(\xiz) 
=\Phi_t\big( f(\xi, z) \ast_t \Phi_t^{-1}(u(\xiz))\big) \\
&=\Phi_t( f(\xi, z))\, \Phi_t\big(\Phi_t^{-1}( u(\xi, z))\big)
 =\Phi_t( f(\xi, z))\, u(\xiz). \nno
\end{align*}
Hence, the proposition follows.
\epfv

By the proposition above and Eqs.\,(\ref{L2.1-e1})  
and (\ref{L2.1-e2}), we also have 
the following corollary. 

\begin{corol}\label{Phi_*-P2-2}
For any $t\in \bC$, $\lambda(\xi)\in \bC[\xi]$ and   
$p(z)\in \cz$, $(\Phi_t)_*$ maps 
the multiplication operators of $\cB_t[\xiz]$   
by $\lambda(\xi)$ and  $p(z)$ 
$($with respect to the product $\ast_t$$)$ 
to the multiplication operators of $\cA[\xiz]$
by $\lambda(\xi)$ and  $p(z)$ 
$($with respect to the product of $\cA[\xiz]$$)$, 
respectively.   
\end{corol}

Note that, as pointed out before, 
the algebras $\cB_t[\xiz]$ $(t\in \bC)$ 
share the same base vectors space $\cxz$. 
Therefore, their Weyl algebras 
$\cD_t[\xiz]$ $(t\in \bC)$ are 
all subalgebras of the algebra 
of linear endomorphisms of $\cxz$. 
The following corollary says that 
all these subalgebras turn  
out to be same, i.e.\@ they do not depend 
on the parameter $t\in \bC$. 

\begin{corol}
For any $t\in \bC$, as subalgebras of the algebra 
of linear endomorphisms of $\cxz$, 
$\cD_t[\xiz]=\cD[\xiz]$. 
\end{corol}
\pf By Proposition \ref{freeness}, $(c)$, 
we know that $\cB_t[\xiz]$ is  
a commuative free algebra generated freely 
by $\xi$ and $z$. By Proposition 
\ref{Phi_*-P1}, $(a)$, 
we know that $\p_i$ 
and $\dlt_i$ $(\lin)$ are 
also derivations 
of $\cB_t[\xiz]$. 
Therefore, the Weyl algebra 
$\cD_t[\xiz]$ as an associative 
algebra over $\bC$ is generated by 
the derivations $\p_i$, $\dlt_i$ $(\lin)$
and the multiplication operators 
(with respect to the product 
$\ast_t$) by $\xi_i, z_i\in \cB_t[\xiz]$ 
$(\lin)$.

By Eqs.\,(\ref{B-Propo-e4}) and (\ref{B-Propo-e5}), 
we see that the multiplication operators 
by $\xi_i, z_i\in \cB_t[\xiz]$ $(\lin)$ 
are same as the operators $\xi_i-t\p_i$ 
and $z_i-t\dlt_i$ which lie in $\cD[\xiz]$.
Hence we have $\cD_t[\xiz] \subseteq \cD[\xiz]$.

To show $\cD[\xiz] \subseteq \cD_t[\xiz]$, 
by Proposition \ref{Phi_*-P1}, $(a)$,
it will be enough 
to show that the multiplication 
operators (with respect 
to the product of $\cA[\xiz]$) 
by $\xi_i, z_i\in \cA[\xiz]$ 
$(\lin)$ also belong to 
$\cD_t[\xiz]$.

But, for any $f(\xiz)\in \bC[\xiz]$, 
by Eqs.\,(\ref{B-Propo-e4}) and (\ref{B-Propo-e5}), 
we have
\begin{align*}
\xi_i f(\xiz)&=(\xi_i-t\p_i) f(\xiz)+t\p_i f(\xiz)
=\xi_i \ast_t f(\xiz)+t\p_i f(\xiz),\\
z_i f(\xiz)&=(z_i-t\dlt_i) f(\xiz)+t\dlt_i f(\xiz)
=z_i \ast_t f(\xiz)+t\dlt_i f(\xiz).
\end{align*}
From the equations above, we see that 
the multiplication operators (with respect 
to the product of $\cA[\xiz]$) 
by $\xi_i, z_i\in \cA[\xiz]$ 
$(\lin)$  do belong  to 
$\cD_t[\xiz]$.
\epfv

\renewcommand{\theequation}{\thesection.\arabic{equation}}
\renewcommand{\therema}{\thesection.\arabic{rema}}
\setcounter{equation}{0}
\setcounter{rema}{0}

\section{\bf  Connections with Interchanges of 
Right and Left Total Symbols of Differential 
Operators} \label{S3}

In this section, we show in Theorem \ref{Interchange} 
that the isomorphisms $\Phi_t$ with $t=\pm 1$ 
coincide with the interchanges between total left and right symbols of differential operators of the polynomial algebra $\cA[z]$.

First, let us fix the following notation and convention for the differential operators of $\cA[z]$.

We denote by $\cD[z]$ the differential operator algebra or the Weyl algebra of $\cA[z]$. For any differential operator $\phi\in \cD[z]$ and polynomial $u(z)\in \cA[z]$, the notation $\phi\, u(z)$ usually denotes the composition of 
$\phi$ and the multiplication operator by $u(z)$. 
So $\phi \, u(z)$ is still a differential operator 
of $\cA[z]$. The polynomial obtained by applying 
$\phi$ to $u(z)$ will be denoted by $\phi(u(z))$. 

Next, let us recall the right and left total 
symbols of differential operators of 
the polynomial algebra $\cA[z]$.
 
For any $\phi \in \cD[z]$, 
it is well-known (e.g. see  
Proposition 2.2 (pp.\,4) in \cite{B} 
or Theorem 3.1 (pp.\,58) in \cite{C}) that 
$\phi$ can be written uniquely as the following 
two finite sums: 
\begin{align}\label{generic-diff}
\phi=\sum_{\alpha \in \bN^n} a_\alpha(z) \p^\alpha
=\sum_{\beta \in \bN^n} \p^\beta b_\beta(z)
\end{align}  
where $a_\alpha(z), b_\beta(z) \in \bC[z]$ but  
denote the multiplication operators by 
$a_\alpha(z)$ and $b_\beta(z)$, respectively.

For the differential operator $\phi\in \cD[z]$ 
in Eq.\,(\ref{generic-diff}), 
the {\it right} and {\it left total symbols} are defined 
to be the polynomials $\sum_{\alpha\in \bN^n} a_\alpha(z) 
\xi^\alpha\in \bC[\xiz]$ and 
$\sum_{\beta \in \bN^n} b_\beta(z) 
\xi^\beta\in \bC[\xiz]$, respectively.
We denote by $\cR: \cD[z] \to \bC[\xi, z]$
(resp.\@ $\cL: \cD[z] \to \bC[\xi, z]$) the linear map 
which maps any $\phi\in \cD[z]$ to its right total symbol 
(resp.\@ left total symbol). 

Note that, by the uniqueness of the expressions in 
Eq.\,(\ref{generic-diff}), both $\cR$ and $\cL$ 
are isomorphisms of vector spaces over $\bC$. 
The interchange of the left (resp.\@ right) total symbol of differential operators  to their right (resp.\@ left) total symbols is given by the isomorphism $\cR \circ \cL^{-1}$ 
(resp.\@ $\cL \circ \cR^{-1}$) from $\cxz$ to $\cxz$.

The main result of this section is the following theorem.

\begin{theo}\label{Interchange} 
As linear maps from $\cxz$ to $\cxz$, we have
\begin{align}
\Phi &=\cR \circ \cL^{-1}.\label{Interchange-e1} \\
\Phi_{t=-1}&=\cL \circ \cR^{-1}. \label{Interchange-e2}
\end{align}
\end{theo}

\pf Note first that, Eq.\,(\ref{Interchange-e2}) follows from Eq.\,(\ref{Interchange-e1}) and the fact that 
$\Phi_{t=-1}=\Phi_{t=1}^{-1}=\Phi^{-1}$.

To show Eq.\,(\ref{Interchange-e1}), since  
both $\Phi$ and $\cR \circ \cL^{-1}$ 
are linear maps, it is enough 
to show that, 
for any $\alpha, \beta\in \bN^n$,  
we have
\begin{align}
\Phi (\xi^\alpha z^\beta)=(\cR \circ \cL^{-1})(\xi^\alpha z^\beta).\label{Interchange-pe1} 
\end{align}

Since 
\begin{align}
(\cR \circ \cL^{-1}) (\xi^\alpha z^\beta)
= \cR( \p^\alpha z^\beta),
\end{align}
so we have to find the right total symbol of the differential operator $\p^\alpha z^\beta\in \cD[z]$.

Note that, for any dummy $u(z)\in \cz$, by the Leibniz rule, we have
\begin{align}
\p^\alpha (z^\beta u(z))&= 
\sum_{\gamma\in \bN^n} 
\binom \alpha\gamma (\p^\gamma z^\beta) 
(\p^{\alpha-\gamma} u(z))\\
&= \left(
\sum_{\gamma\in \bN^n} 
\binom{\alpha}{\gamma} (\p^\gamma z^\beta) 
\p^{\alpha-\gamma} \right) u(z).\nno
\end{align}
Therefore, the right total symbol of the differential operator 
$\p^\alpha z^\beta \in \cD[z]$ is given by 
\begin{align*}
\cR (\p^\alpha z^\beta) &= 
\sum_{\gamma\in \bN^n} 
\binom{\alpha}{\gamma} (\p^\gamma z^\beta) 
\xi^{\alpha-\gamma} =
\sum_{\gamma\in \bN^n}
\binom{\alpha}{\gamma}  \xi^{\alpha-\gamma} (\p^\gamma z^\beta)  \\
&=\sum_{\gamma\in \bN^n} 
\frac 1{\gamma!} (\dlt^\gamma \xi^\alpha)(\p^\gamma z^\beta) \nno
\end{align*}

Combining the equation above with 
Eqs.\,(\ref{B-Propo-e3}) and (\ref{L2.1-e3}) with $t=-1$, 
we have 
\begin{align}
\cR (\p^\alpha z^\beta)
= \xi^\alpha \ast_{t=-1} z^\beta 
= \Phi_{t=1}( \xi^\alpha z^\beta)= 
\Phi( \xi^\alpha z^\beta).
\end{align}

Hence, we have proved Eq.\,(\ref{Interchange-pe1}) 
and also the theorem.
\epfv

\begin{corol}\label{C3.2}
For any $\lambda(\xi) \in \bC[\xi]$ and $p(z) \in \bC[z]$, 
we have
\begin{align}
\cR ( \lambda(\p) p(z))&= \lambda(\xi) \ast_{t=-1} p(z). 
\label{T-R2L}\\
\cL (p(z)\lambda(\p)) &= \lambda(\xi) \ast p(z). 
\label{T-L2R}
\end{align}
\end{corol}

\pf By Eqs.\,(\ref{Interchange-e1}) and 
(\ref{L2.1-e3}) with $t=1$, we have 
\begin{align*}
\cR ( \lambda(\p) p(z))&= \cR (\cL^{-1} (\lambda(\xi) p(z)))  
=(\cR\circ \cL^{-1}) (\lambda(\xi) p(z)) \\
&= \Phi_{t=1} (\lambda(\xi) p(z)) 
=\lambda(\xi) \ast_{t=-1} p(z).
\end{align*}

So we have Eq.\,(\ref{T-R2L}). 
Eq.\,(\ref{T-L2R}) can be proved similarly
by using Eqs.\,(\ref{Interchange-e2}) and 
(\ref{L2.1-e3}) with $t=-1$.
\epfv

Finally, we end this section with the following 
one-variable example.

\begin{exam}
Let $n=1$ and $\phi= z^2 \p^3$. Then, 
\begin{align*}
\cR(\phi)& = \cR( z^2 \p^3)=\xi^3 z^2. \\
\cL(\phi)&=\cL(z^2\p^3)
= \xi^3 \ast z^2 
=(z-\dlt)^2 \xi^3 \\
&=(z^2-2z\dlt+\dlt^2) \xi^3 
=\xi^3z^2- 6\xi^2 z+6\xi. 
\end{align*}
Therefore, we have 
\begin{align*}
\phi = z^2 \p^3=\p^3 z^2- 6\p^2 z+6\p.
\end{align*}
\end{exam}

\renewcommand{\theequation}{\thesection.\arabic{equation}}
\renewcommand{\therema}{\thesection.\arabic{rema}}
\setcounter{equation}{0}
\setcounter{rema}{0}

\section{\bf  A Re-formulation of the Image Conjecture on Commuting Differential Operators of Order One with Constant Leading Coefficients}\label{S4}

In this section, we show that the algebra $\cB_t[\xiz]$ 
$(t\in \bC)$ is closely related with a theorem (See 
Theorem \ref{D_t-Taylor}) first proved in \cite{IC} and also with the so-called {\it image} conjecture 
(See Conjecture \ref{IC1}) proposed in \cite{IC} on the differential operators $\xi-t\p$ $(t\in \bC)$. 

In Subsection \ref{S4.1}, we use certain Taylor series expansion of elements of $\cB_t[\xiz]$ to give a new and more conceptual  proof for Theorem \ref{D_t-Taylor}. In Subsection \ref{S4.2}, we first give a new formulation (See Conjecture \ref{IC2}) for 
Conjecture \ref{IC1} in terms of the algebra 
$\cB_t[\xiz]$ and the notion of {\it Mathieu subspaces} 
(see Definition \ref{M-ideal}) introduced in 
\cite{GIC}, and then show in Theorem \ref{IC1=IC2} 
that the new formulation is indeed equivalent 
to Conjecture \ref{IC1}.

\subsection{The Taylor Series with Respect to the Product $\ast_t$}
\label{S4.1}

First, let us recall the following elementary fact on polynomials in $\xi$ and $z$.

For any $f(\xiz) \in \cA[\xiz]$, 
we may view $f(\xiz)$ as a polynomial in $\xi$ 
with coefficients in $\cA[z]$. Then it has 
the following Taylor series expansion 
\begin{align}\label{Taylor-e1}
f(\xiz)=\sum_{\alpha\in \bN^n} 
\frac 1{\alpha!}\,\, \xi^\alpha  c_\alpha(z)
\end{align}
for some $c_\alpha(z)\in \cA[z]$.

Let $ev_{{}_0}: \cA[\xiz]\to \cA[z]$ be 
the {\it evaluation} map of $\cA[\xiz]$ 
at $\xi=0$, i.e.\@ for any 
$u(\xiz)\in \cA[\xiz]$, 
$ev_{{}_0}(u)\!:=u(0, z)$. Then, the $c_\alpha(z)$ 
$(\alpha\in \bN^n)$ in Eq.\,(\ref{Taylor-e1}) 
are given by
\begin{align}\label{for4c-alpha}
c_\alpha(z)=ev_{{}_0}(\delta^\alpha f).
\end{align}

Note that another characterization of 
the evaluation map $ev_{{}_0}$ is that 
$ev_{{}_0}$ is the (unique) algebra homomorphism 
from $\cA[\xiz]$ to $\cA[z]$ with $ev_{{}_0}(\xi_i)=0$ 
and $ev_{{}_0}(z_i)=z_i$ 
for any $1\le i\le n$.

Now, come back to our algebras  
$\cB_t[\xiz]$ $(t\in \bC)$. 
By Proposition \ref{freeness}, 
$(c)$, we know that 
$\cB_t[\xiz]$ is also 
a commutative free algebra 
generated freely by $\xi$ and $z$ 
with the same base 
vector space $\cxz$. 
Hence, we should expect 
similar expansions as in 
Eq.\,(\ref{Taylor-e1}) 
for polynomials $f(\xiz)\in \cxz$ 
with respect to the product $\ast_t$.

But, in order to formulate the expected 
expansions precisely, we need first to 
introduce the analogue of 
the evaluation map $ev_{{}_0}$ 
for the algebra $\cB_t[\xiz]$.

Note that, by Proposition \ref{freeness}, 
$(b)$, the subalgebra of $\cB_t[\xiz]$ 
generated by $z$ is also 
$\cA[z]\subset \bC[\xiz]$. 
Parallel to the second characterization 
of the evaluation map $ev_{{}_0}$ 
mentioned above,  we let $\cE_t$ 
be the unique algebra 
homomorphism from $\cB_t[\xiz]\to \cA[z]$ 
such that $\cE_t(\xi_i)=0$ and 
$\cE_t(z_i)=z_i$ for any $\lin$.

Note also that, by   
Eqs.\,(\ref{L2.1-e1}) and (\ref{L2.1-e2}), 
the algebra isomorphism
$\Phi_t: \cB_t[\xiz]\to \cA[\xiz]$ 
maps $\xi_i$ (resp.\@ $z_i$) to  
$\xi_i$ (resp.\@ $z_i$) for 
any $\lin$. Hence the composition 
$ev_{{}_0}\circ \Phi_t:\cB_t[\xiz]\to \cA[z]$ 
has the same characterizing property of $\cE_t$. 
Therefore, we have 
\begin{align}\label{cE=ev}
\cE_t=ev_{{}_0}\circ \Phi_t.
\end{align} 

Furthermore, we can also derive a 
more explicit formula for $\cE_t$ as follows. 

For any $\alpha\in \bN^n$ and 
$p(z)\in \bC[z]$, consider
\begin{align}\label{ExplicitEt}
\cE_t(\xi^\alpha p(z))
&=ev_{{}_0} ( \Phi_t( \xi^\alpha p(z) ) ) \\
\intertext{Applying Eq.\,(\ref{L2.1-e3}) and then Eq.\,(\ref{B-Propo-e4}) with $t$ replaced by $-t$:}
&=ev_{{}_0} ( \xi^\alpha \ast_{-t} p(z) ) ) \nno \\
&=ev_{{}_0} ( (\xi+t\p)^\alpha ( p(z)) ) 
=t^{|\alpha|} \p^\alpha (p(z)). \nno
\end{align}

From the formula above, we see that,  
for any $g(z, \xi)\in \bC[z, \xi]$, 
$\cE_t (g(z, \xi))\in \bC[z]$ can be obtained by, 
first, writing each monomial of $g(z, \xi)$ as  
$\xi^\beta z^\gamma$ $(\beta, \gamma\in \bN^n)$, 
i.e.\@ putting the free variables $\xi_i$'s to the most left 
in each monomial of $g(z, \xi)$, 
and then replacing the part $\xi^\beta$ by 
the differential operator $t^{|\beta|} \p^\beta$ and 
applying it to the other part
$z^\gamma$ of the monomial. 
For examples, we have 
{\it 
\begin{align*}
\cE_t(1)&=1; \\
\cE_t(z^\alpha)&=(t\p)^0(z^\alpha)=z^\alpha \qquad \quad  \, \, \qquad\qquad\mbox{for any } \alpha\in \bN^n; \\
\cE_t(\xi^\alpha)&= t^{|\alpha|} \p^\alpha(1)=0  \qquad\qquad \qquad \quad \quad  \mbox{for any } 0\neq\alpha\in \bN^n; \\
\cE_t(z_1^m \xi_1^2)&=t^2\p_1^2(z_1^m)=m(m-1)t^2 z_1^{m-2}  \qquad \, \mbox{for any } m\ge 2.  
\end{align*}  }

Now we are ready to formulate and prove 
the expected expansion of polynomials with respect to the new product $\ast_t$, which is parallel to the Taylor expansion in Eq.\,(\ref{Taylor-e1}). 

\begin{theo}\label{D_t-Taylor}
For any $t\in \bC$ and $f(\xi, z)\in \cxz$, we have 
\begin{align}
f(\xi, z)&=\sum_{\alpha\in \bN^n} \frac 1{\alpha!} \,\,
\xi^\alpha \ast_t a_\alpha(z),\label{D_t-Taylor-e1} \\
f(\xi, z)&=\sum_{\alpha\in \bN^n} \frac 1{\alpha!} \,
(\xi-t\p_z)^\alpha a_\alpha(z),\label{D_t-Taylor-e2}
\end{align}
where, for any $\alpha\in \bN^n$,
\begin{align}\label{D_t-Taylor-e3}
a_\alpha(z)=\cE_t(\dlt^\alpha f).
\end{align}
Furthermore, the expansions of the forms  
in Eqs.\,$(\ref{D_t-Taylor-e1})$
and $(\ref{D_t-Taylor-e2})$ for $f(\xiz)$ 
are unique.
\end{theo}

\pf Note first that, 
by Eq.\,(\ref{B-Propo-e4}) 
in Proposition \ref{B-Propo}, 
Eq.\,(\ref{D_t-Taylor-e1})
and Eq.\,(\ref{D_t-Taylor-e2}) 
are actually equivalent. 
So we will focus only on 
Eq.\,(\ref{D_t-Taylor-e1}).

The uniqueness of the expansion in 
Eq.\,(\ref{D_t-Taylor-e1}) follows 
directly from Proposition \ref{freeness}, 
$(a)$-$(c)$. 

To show that Eq.\,(\ref{D_t-Taylor-e1}) with
$a_\alpha(z)$ $(\alpha\in \bN^n)$ 
given in Eq.\,(\ref{D_t-Taylor-e3})
does hold, we first write the 
Taylor series expansion 
of $\Phi_t(f(\xiz))$ as in Eq.\,(\ref{Taylor-e1}): 
\begin{align}\label{D_t-Taylor-pe2}
\Phi_t (f(\xiz))=\sum_{\alpha\in \bN^n} 
\frac 1{\alpha!}\, \, \xi^\alpha a_\alpha(z)
\end{align}
where $a_\alpha(z)\in \cz$ $(\alpha\in \bN^n)$ 
are given by 
\begin{align}\label{D_t-Taylor-pe3}
a_\alpha(z)=ev_{{}_0}(\dlt^\alpha \Phi_t (f)).
\end{align}

Applying $\Phi_{-t}$ to Eq.(\ref{D_t-Taylor-pe2}) and, 
by Eq.\,(\ref{L2.1-e3}) with $t$ replaced by $-t$, 
we get Eq.\,(\ref{D_t-Taylor-e1}). 

Next, note that $\delta^\alpha$ $(\alpha\in \bN^n)$ 
commute with $\Lambda$, hence they also 
commute with $\Phi_t=e^{t\Lambda}$. Then, 
by Eqs.\,(\ref{D_t-Taylor-pe3}) and 
(\ref{cE=ev}), we have 
\begin{align*}
a_\alpha(z)=ev_{{}_0}(\Phi_t (\dlt^\alpha f) )
=(ev_{{}_0}\circ\Phi_t) (\dlt^\alpha f)=\cE_t(\dlt^\alpha f).
\end{align*}
Therefore, Eq.\,(\ref{D_t-Taylor-e3}) also holds.
\epfv

Several remarks on Theorem \ref{D_t-Taylor} 
and the proof above are as follows.

First, Theorem \ref{D_t-Taylor} with $t=1$ 
was first proved in \cite{IC}.
The proof in \cite{IC} is more straightforward. 
It does not use the algebra $\cB_t[\xiz]$ 
and the product $\ast_t$. But the proof given 
here is more conceptual. For example, 
the expansion in Eq.\,$(\ref{D_t-Taylor-e2})$ 
becomes much more natural after we show here that 
it is just the usual Taylor series expansion 
of polynomials as in Eq.\,$(\ref{Taylor-e1})$ but 
in the new context of the algebra $\cB_t[\xiz]$.

Second, Eq.\,(\ref{D_t-Taylor-e3}) can also be derived 
directly from  Eq.\,(\ref{D_t-Taylor-e2}) as 
in \cite{IC}. Namely, apply $\dlt^\alpha$ 
to Eq.\,(\ref{D_t-Taylor-e2}) and then replace 
$\xi$ by $t\p$ in both sides 
of the resulting equation.

Third, not all formal power series $f(\xiz)\in \cA[[\xiz]]$ 
can be expanded in the form of Eq.\,$(\ref{D_t-Taylor-e1})$
or $(\ref{D_t-Taylor-e2})$. For example, let $n=1$ and 
$f(\xiz)=e^{\xi z}$ and assume  
that $(\ref{D_t-Taylor-e2})$ holds 
for $f(\xiz)$. Then, by the argument 
in the previous paragraph, we see that 
$a_m(z)$ $(m\ge 0)$ must be given by   
Eq.\,(\ref{D_t-Taylor-e3}). 
But, for the series $\dlt^m f(\xiz)=z^m 
\sum_{k\ge 0} \frac {(\xi z)^k}{k!}$, 
$\cE_t$ is not well-defined, which is 
a contradiction.

Another way to look at the fact 
above is as follows. 
Even though $\cB_t[\xiz]$ $(t\ne 0)$ 
and $\cA[\xiz]$ share the same 
base vector space $\cxz$, 
by Proposition \ref{L-topology}, 
we know that they induce different 
$\ell$-adic topologies 
on $\cxz$. Therefore, 
their completions with respect 
to the different $\ell$-adic 
topologies will be 
different. In other words, 
the formal power series algebras  
with respect to the product 
$\ast_t$ $(t\ne 0)$ 
and the usual formal power series 
algebra $\cA[[\xiz]]$ do not 
share the same base vector 
space anymore. 

For the example $f(\xiz)=e^{\xi z}$ 
above, we have $f(\xiz)\in \cA[[\xiz]]$. 
But, by the argument in the proof of 
Proposition \ref{L-topology} with 
$\alpha_m$ $(m\ge 1)$ replaced by $m$, 
it is easy to see that, for any $t\ne 0$, 
$f(\xiz)=e^{\xi z}$ does not 
lie in the completion of 
$\cB_t[\xiz]$ with respect 
to the $\ell$-adic 
topology on $\bC[\xiz]$ induced by 
$\cB_t[\xiz]$. Therefore, 
$f(\xiz)=e^{\xi z}$ can not be written 
as a formal power series with respect to 
the product $\ast_t$ as in 
Eq.\,(\ref{D_t-Taylor-e1}).

\subsection{Re-Formulation of the Image Conjecture in Terms of 
the Algebra $\cB_t[\xiz]$}\label{S4.2}

First let us recall the following notion 
introduced recently in \cite{GIC}.

\begin{defi}\label{M-ideal}
Let $R$ be any commutative ring and $\cA$ 
a commutative $R$-algebra.
We say that an $R$-subspace 
$\cM$ of $\cA$ is a {\it Mathieu subspace} of $\cA$  
if the following property holds: 
if $a\in \cA$ satisfies $a^m\in \cM$ for all $m\ge 1$,
then, for any $b\in\cA$, we have 
$ b a^m \in \cM$ for all $m\gg 0$, 
i.e.\@ there exists $N\ge 1$ $($depending on $a$ and $b$$)$ 
such that $b a^m  \in \cM$ for all $m\ge N$.
\end{defi}   

From the definition above, it is easy to see that any ideal of 
$\cA$ is automatically a Mathieu subspace of $\cA$, 
but not conversely (See \cite{GIC} for some examples of Mathieu subspaces 
which are not ideals). Therefore, the notion of Mathieu subspaces 
can be viewed as a generalization of the notion of ideals.

Next, for any $t\in \bC$, set 
\begin{align}\label{Def-im}
\im(\xi-t\p)\!:=\sum_{i=1}^n (\xi_i-t\p_i)\bC[\xiz].
\end{align}

We call $\im(\xi-t\p)$ the {\it image} of the 
commuting differential operators 
$(\xi_i-t\p_i)$ $(\lin)$. 

With the notion and notation fixed above, 
the {\it image} conjecture proposed in \cite{GIC} 
for the commuting differential operators $(\xi-t\p)$ 
can be re-stated as follows.

\begin{conj}\label{IC1}
For any $t\in \bC$, $\im(\xi-t\p)$ is 
a Mathieu subspace of the polynomial 
algebra $\cA[\xiz]$.
\end{conj}

One of the motivations of the conjecture above is 
the following theorem proved in \cite{IC}.
 
\begin{theo}\label{IC1=>JC}
Conjecture \ref{IC1} implies 
the Jacobian conjecture.
\end{theo}

Actually, it has been shown in \cite{IC} that the {\it Jacobian} 
conjecture is equivalent to some very special cases of 
Conjecture \ref{IC1}. For more detail, see \cite{IC}.

The main result of this subsection is to 
show that the conjecture above can actually 
be re-formulated as follows.

\begin{conj}\label{IC2}
Set $\xi\bC[\xiz]\!:=\sum_{i=1}^m \xi_i\cxz$. Then, for any $t\in \bC$, $\xi\bC[\xiz]$ as a subspace of 
$\cB_t[\xiz]$ is a Mathieu subspace of $\cB_t[\xiz]$.
\end{conj}

\begin{theo}\label{IC1=IC2}
Conjecture \ref{IC2} is equivalent 
to Conjecture \ref{IC1}.
\end{theo}

\pf First, denote by $\xi \ast_t\bC[\xiz]$ the ideal of 
$\cB_t[\xiz]$ generated by $\xi_i$ $(\lin)$. 
View $\xi \ast_t\bC[\xiz]$ as a subspace of $\cA[\xiz]$ 
and apply Eqs.\,(\ref{Def-im}) and (\ref{B-Propo-e4}), 
we have  
\begin{align}\label{IC1=IC2-pe1}
\im(\xi-t\p)=\sum_{i=1}^n \xi_i \ast_t \bC[\xiz]
=\xi \ast_t\bC[\xiz].
\end{align}

Second, by Eqs.\,(\ref{GD-trivial-e1}) 
and (\ref{L2.1-e1}), we have
\begin{align*}
\Phi_t( \xi \ast_t \bC[\xiz])= 
\Phi_t(\xi)\Phi_t( \bC[\xiz]) 
=\xi \bC[\xiz]
\end{align*}

Hence, we also have 
\begin{align}\label{IC1=IC2-pe2}
\xi \ast_t \bC[\xiz]
=\Phi_t^{-1}(\xi \bC[\xiz])=\Phi_{-t}(\xi \bC[\xiz]).
\end{align}

Combine Eqs.\,(\ref{IC1=IC2-pe1}) 
and (\ref{IC1=IC2-pe2}), we get 
\begin{align}\label{IC1=IC2-pe3}
\Phi_{-t}(\xi \bC[\xiz])=\im(\xi-t\p).
\end{align}
 
Third, by Proposition $4.9$ in \cite{GIC}, 
we know that pre-images of 
Mathieu subspaces under algebra 
homomorphisms are still 
Mathieu subspaces, from which    
it is easy to check that  
Mathieu subspaces are  
preserved by algebra 
isomorphisms.  
By using this fact 
(on the algebra isomorphism 
$\Phi_{-t}:\cB_{-t}[\xiz]\to \cA[\xiz]$) 
and also Eq.\,(\ref{IC1=IC2-pe3}), we see that,  
$\xi \bC[\xiz]$ is a Mathieu subspace of 
$\cB_{-t}[\xiz]$ iff 
$\im(\xi-t\p)$ is a Mathieu subspace 
of $\cA[\xiz]$.

Replacing $t$ by $-t$ in the equivalence above, 
we have that, $\xi \bC[\xiz]$ is a Mathieu subspace of 
$\cB_t[\xiz]$ for any $t\in \bC$ iff $\im(\xi+t\p)$ is a Mathieu subspace of $\cA[\xiz]$ for any $t\in \bC$ 
iff $\im(\xi-t\p)$ 
is a Mathieu subspace of $\cA[\xiz]$ for any $t\in \bC$. Hence, we have proved the theorem. 
\epfv

From Theorems \ref{IC1=>JC} and \ref{IC1=IC2}, 
we immediately have the following corollary.

\begin{corol}\label{IC2=>JC}
Conjecture \ref{IC2} implies 
the Jacobian conjecture.
\end{corol}

\begin{rmk}\label{rmk-S4.2-1}
Note that, when $t=0$, Conjecture \ref{IC2} 
is trivial since $\xi\cxz$ is an ideal of 
the algebra $\cB_{t=0}[\xiz]=\cA[\xiz]$. 
In general, Conjecture \ref{IC2} in some 
sense just claims that the algebras  
$\cB_t[\xiz]$ $(t\in \bC)$ do not differ or change 
too much from $\cA[\xiz]$ so that 
the vector subspace $\xi\cxz$ still 
remains as a Mathieu subspace 
of $\cB_t[\xiz]$. 

From this point of view, 
the triviality of 
the deformation 
$\cB_t[\xiz]$ $(t\in \bC)$
of the polynomial algebra $\cA[\xiz]$ 
given in Corollary  
\ref{D-trivial} may be viewed 
as a fact in favor of   
Conjecture \ref{IC2}, 
hence also to the Jacobian conjecture 
via the implication 
in Corollary \ref{IC2=>JC}.
\end{rmk}

\begin{rmk}\label{rmk-S4.2-2}
Conjecture \ref{IC2} and also the Jacobian conjecture can be viewed as problems caused by the following fact. 
Namely, due to the change  of the algebra 
structure from $\cA[\xiz]$ 
to $\cB_t[\xiz]$, the evaluation map at $\xi=0$, 
which is an algebra homomorphism from $\cA[\xiz]$ to $\cA[z]$, is not an algebra homomorphism from $\cB_t[\xiz]$ to $\cA[z]$ if $t\ne 0$. Therefore, its kernel $\xi\bC[\xiz]$ does not remain to be 
an ideal of $\cB_t[\xiz]$ anymore. 

But, on the other hand, as we will see later in 
Subsection \ref{S5.2} $($See Theorem \ref{ast2LP} 
and Remark \ref{rmk-S5.2-1}$)$, the same fact 
for the evaluation map at $\xi=1$, i.e.\@ $\xi_i=1$ 
$(\lin)$, in some sense also causes something truely  remarkable, namely, the generalized 
Laguerre polynomials.
\end{rmk}

\renewcommand{\theequation}{\thesection.\arabic{equation}}
\renewcommand{\therema}{\thesection.\arabic{rema}}
\setcounter{equation}{0}
\setcounter{rema}{0}

\section{\bf  Connections with the Generalized Laguerre Polynomials}
\label{S5}

In this section, we study some connections and interactions of the monomials of the algebra $\cB[\xiz]$ in $\xi$ and $z$ with the generalized Laguerre polynomials in one or more variables.

In Subsection \ref{S5.1}, we briefly recall the definition 
and the orthogonal property of the generalized 
Laguerre polynomials. In Subsection \ref{S5.2}, 
we show that the generalized Laguerre polynomials 
can be obtained from  
certain monomials of the algebra $\cB[\xiz]$ 
in $\xi$ and $z$ (See Theorem \ref{ast2LP} and 
Corollary \ref{ast2LP-corol}).  

In Subsection \ref{S5.3}, we study some 
applications of the connection given in 
Theorem \ref{ast2LP}.  We first use certain 
properties of the generalized 
Laguerre polynomials to derive some results on   
some monomials of $\cB[\xiz]$ in $\xi$ and $z$.
We then use some results derived in Section \ref{S2} 
on the algebra $\cB[\xiz]$ to give new proofs for 
some important properties of the generalized 
Laguerre polynomials (see Proposition \ref{Recur-Formula} 
and Theorem \ref{Aso-Lag-Equation}).

\subsection{The Generalized Laguerre Orthogonal Polynomials} 
\label{S5.1}

First, let us recall the generalized Laguerre 
orthogonal polynomials in one variable.

For any $k\in \bR$ and $m\in \bN$, the generalized Laguerre polynomial $L_m^{[k]}(z)$ in one variable is given by 
\begin{align}\label{Explicit-Formula}
L_m^{[k]}(z) =\sum_{j=0}^m \binom{m+k}{m-j}\, 
\frac{(-z)^j}{j!}.
\end{align}

Here we are only interested in the case that $k\in \bN$. 
For any fixed $k\in \bN$, the generating function of  
the generalized Laguerre polynomials $L_m^{[k]}(z)$ $(m\ge 0)$ is given by
\begin{align}\label{Generatn-Fctn}
\frac{ \exp(-\frac{z u}{1-u})}{ (1-u)^{k+1} }
=\sum_{m=0}^{+\infty} L_m^{[k]}(z) \, u^m,
\end{align}
where $u$ above denotes a formal variable 
which commutes with $z$.

The multi-variable generalized Laguerre 
polynomials are defined as follows.

Let ${\bf k}=(k_1, k_2, ..., k_n)\in \bN^n$ and 
${\alpha}=(a_1, a_2, ..., a_n)\in \bN^n$. 
The generalized Laguerre polynomials in $n$-variable 
$z=(z_1, z_2, ..., z_n)$ is 
defined by 
\begin{align}\label{Mult-LPs}
L_\alpha^{[{\bf k}]}(z)\!:=L_{a_1}^{[k_1]}(z_1)
L_{a_2}^{[k_2]}(z_2)\cdots L_{a_n}^{[k_n]}(z_n).
\end{align}

The polynomials $L_\alpha(z)\!:=L_\alpha^{[0]}(z)$ 
$(\alpha\in \bN^n)$ are the so-called the (classical) 
{\it Laguerre} polynomials. They were named after 
Edmond. N. Laguerre \cite{L}. The generalized 
Laguerre polynomials were introduced much later 
by G. P\'olya and G. Szeg\"o \cite{PS} in $1976$. 

One of the most important properties of the generalized Laguerre polynomials is the following theorem.

\begin{theo}\label{ortho}
For any ${\bf k}, \alpha, \beta \in \bN^n$, we have 
\begin{align}\label{ortho-e1}
\int_{(\bR_{>0})^n} 
L_\alpha^{[{\bf k}]} (z) L_\beta^{[{\bf k}]} (z) 
\, w(z)\, dz 
=\delta_{\alpha, \beta}\, \frac{(\alpha+{\bf k})!}{\alpha!},
\end{align}
where $\delta_{\alpha, \beta}$ is the 
Kronecker delta function and 
$w(z)$ given by 
\begin{align}\label{Lgrwt}
w(z)\!:=z^{\bf k} e^{-\sum_{i=1}^n z_i}.
\end{align}
\end{theo}

The function $w(z)$ above is called the {\it weight} 
function of the generalized Laguerre polynomials 
$L_\alpha^{[\bf k]}(z)$ $(\alpha\in \bN^n)$.

Consequently, with any fixed ${\bf k}$, 
the generalized Laguerre polynomials 
$L_\alpha^{[{\bf k}]} (z)$ $(\alpha \in \bN^n)$
form an orthogonal basis of $\cz$ with 
respect to the Hermitian form defined by 
\begin{align}
(f, g)=\int_{(\bR_{>0})^n} 
f(z)\bar g(z) \, w(z)\, dz, 
\end{align}
where $\bar g(z)$ denotes the complex conjugation of the polynomial $g(z)\in \cz$.

There are many other interesting and 
important properties of the generalized 
Laguerre polynomials. We refer the reader 
to \cite{Sz}, \cite{PS}, \cite{AAR} 
and \cite{DX} for very thorough study on 
this family of orthogonal polynomials.
See also the {\it Wolfram Research}  
web sources \cite{W1} and \cite{W2} for 
over one hundred formulas and identities  
on the (generalized) Laguerre polynomials.

\subsection{The Generalized Laguerre Polynomials 
in Terms of the Product $\ast$}\label{S5.2}

The main result of this subsection is the following theorem.

\begin{theo}\label{ast2LP}
For any ${\bf k}, \alpha \in \bN^n$, we have 
\begin{align}
L_{\alpha}^{[{\bf k}]}( \xi z)
&=\frac {(-1)^{|\alpha|}}{\alpha!} \, 
 \xi^{-{\bf k}}(\xi^{\alpha+{\bf k}} \ast z^\alpha). 
\label{ast2LP-e1} \\
L_{\alpha}^{[{\bf k}]}( \xi z)&=\frac {(-1)^{|\alpha|}}{\alpha!} \, 
z^{-{\bf k}}(\xi^\alpha \ast z^{\alpha+{\bf k}}), \label{ast2LP-e2}
\end{align}
where $\xi z\!:=(\xi_1z_1, \xi_2z_2, ..., \xi_nz_n)$.

In particular, for the Laguerre polynomials, we have 
\begin{align}\label{ast2LP-e3}
L_\alpha( \xi z)=\frac {(-1)^{|\alpha|}}{\alpha!}\, \, \xi^\alpha \ast z^\alpha.
\end{align}
\end{theo}

\pf We first prove Eq.\,(\ref{ast2LP-e3}). 
Note first that, as pointed out in 
Subsection $2.1$ \cite{GIC},
the Laguerre polynomials $L_m(z)$ 
$(m\in \bN)$ in one variable can be obtained as
\begin{align}
L_m(z)=\frac 1{m!}(\p-1)^m(z^m).
\end{align}

Changing the variable $z\to \xi z$ in the equation above, 
we get
\begin{align*}
L_m(\xi z)&=\frac 1{m!}(\xi^{-1}\p-1)^m(\xi^m z^m)
=\frac 1{m!}\xi^{-m}(\p-\xi)^m(\xi^m z^m)\\
&=\frac 1{m!}(\p-\xi)^m (z^m) 
=\frac {(-1)^m}{m!}(\xi-\p)^m(z^m). \nno 
\end{align*}

By Eq.\,(\ref{Mult-LPs}) with ${\bf k}=0$ 
and the equation above, we see that 
the multi-variable Laguerre polynomials $L_\alpha(z)$ 
$(\alpha\in \bN^n)$ can be given by
\begin{align}
L_\alpha(\xi z)
=\frac {(-1)^{|\alpha|}}{\alpha!}(\xi-\p)^\alpha(z^\alpha).  
\end{align}
Then, apply Eq.\,(\ref{B-Propo-e4}) with 
$\lambda(\xi)=\xi^\alpha$ and 
$t=1$,  we get Eq.\,(\ref{ast2LP-e3}).

To show  Eq.\,(\ref{ast2LP-e1}), recall that we have the following well-known identity for the one-variable generalized Laguerre polynomials, which can be easily derived from the generating functions of the generalized Laguerre polynomials in 
Eq.\,(\ref{Generatn-Fctn}): 
\begin{align}\label{02k}
L_m^{[k]}(z)=(-1)^k \p^k L_{m+k}(z).
\end{align}

Now, by Eq.\,(\ref{Mult-LPs}) and the equation above, we see that 
the multi-variable generalized Laguerre polynomials can be given by
\begin{align}
L_\alpha^{[{\bf k}]}(z)=(-1)^{|\bf k|} 
\p^{\bf k} L_{\alpha+{\bf k}}(z).
\end{align}

Changing the variable $z\to \xi z$ in the equation above, 
we get
\begin{align}
L_\alpha^{[{\bf k}]}(\xi z)&=(-1)^{|\bf k|} 
(\p^{\bf k} L_{\alpha+{\bf k}}) (\xi z) \\
&=(-1)^{|\bf k|}\xi^{-{\bf k}} 
\p^{\bf k} (L_{\alpha+{\bf k}} (\xi z) ) \nno \\
\intertext{Applying Eq.\,(\ref{ast2LP-e3}) and then   
Eq.(\ref{ast-Leibniz-e1}):}
&=\frac{(-1)^{|\alpha|}}{(\alpha+{\bf k})!}
\xi^{-{\bf k}} \p^{\bf k} (\xi^{\alpha+{\bf k}} \ast 
z^{\alpha+{\bf k}}). \nno \\
&=\frac{(-1)^{|\alpha|}}{\alpha!}
\xi^{-{\bf k}} (\xi^{\alpha+{\bf k}} \ast 
z^{\alpha}). \nno
\end{align}

Hence, we get Eq.\,(\ref{ast2LP-e1}). Switching $\xi$ and $z$ in 
Eq.\,(\ref{ast2LP-e1}) and using the commutativity 
of the product $\ast$, we get Eq.\,(\ref{ast2LP-e2}). 
\epfv

\begin{corol}\label{ast2LP-corol}
For any ${\bf k}, \alpha \in \bN^n$, we have 
\begin{align*}
L_\alpha(z)&=\frac {(-1)^{|\alpha|}}{\alpha!}\, \, \left. 
(\xi^\alpha \ast z^\alpha)\right|_{\xi=1}; \\
\\
L_\alpha^{[{\bf k}]}(z)&=\frac {(-1)^{|\alpha|}}{\alpha!} \, 
\left. (\xi^{\alpha+{\bf k}} \ast z^\alpha)\right|_{\xi=1}; \\
\\ 
L_{\alpha}^{[{\bf k}]}(z)&=\frac {(-1)^{|\alpha|}}{\alpha!} \, 
\left. z^{-{\bf k}}(\xi^{\alpha} \ast z^{\alpha+{\bf k}})
\right|_{\xi=1},  
\end{align*}
where $|_{{}_{\xi=1}}$ denotes the evaluation map from 
$\bC[\xiz]$ to $\cz$ by setting $\xi_i=1$ 
for any $\lin$.
\end{corol}

\begin{rmk}\label{rmk-S5.2-1}
Note that, the evaluation map $|_{{}_{\xi=1}}$ 
viewed as a linear map from $\cA[\xiz]$ to $\cA[z]$ 
is a homomorphism of algebras. But, as a linear map from the algebra $\cB[\xiz]$ to the polynomial algebra 
$\cA[z]$, it is not a homomorphism of algebras anymore. 
In particular, we have 
\begin{align*}
\left. (\xi^\alpha \ast z^\alpha)\right|_{\xi=1} 
\ne 1 \ast z^\alpha=z^\alpha.
\end{align*}
Otherwise the generalized Laguerre polynomials would be trivialized. 

Therefore, in some sense,  the fact that 
the evaluation map 
$|_{{}_{\xi=1}}: \cB[\xiz] \to \cA[z]$ 
fails to be an algebra homomorphism causes 
the non-trivial, actually truly remarkable, 
generalized Laguerre polynomials.  
But, on the other hand, as we have discussed in 
Subsection \ref{S4.2} 
$($See Remark \ref{rmk-S4.2-2}$)$, 
the same fact for the evaluation map at $\xi=0$ also 
causes some extremely difficult open problems 
such as Conjecture \ref{IC2} and the Jacobian conjecture.
\end{rmk}

Another immediate consequence of Theorem \ref{ast2LP} 
is the following corollary.

\begin{corol}\label{even}
For any $\alpha, \beta\in \bN^n$, we have
\begin{align}\label{even-e1}
\xi^\beta(\xi^\alpha\ast z^{\alpha+\beta})
=z^\beta(\xi^{\alpha+\beta}\ast z^\alpha).
\end{align}
\end{corol}

Note that the corollary follows immediately from 
Eqs.\,(\ref{ast2LP-e1}) and (\ref{ast2LP-e2})
with ${\bf k}=\beta$. 
But here we also give a more 
straightforward proof.

\pf Consider 
\allowdisplaybreaks{
\begin{align*}
\xi^\beta(\xi^\alpha \ast z^{\alpha+\beta})
&=(\xi-\p+\p)^\beta (\xi^\alpha \ast z^{\alpha+\beta}) \\
&=\sum_{\gamma\in \bN^n} \binom \beta \gamma (\xi-\p)^{\beta-\gamma}\p^\gamma 
(\xi^\alpha \ast z^{\alpha+\beta}) \\
\intertext{Applying Eq.\,(\ref{ast-Leibniz-e1}) and then   
Eq.(\ref{B-Propo-e4}):}
&=\sum_{\gamma\in \bN^n} \binom \beta \gamma 
\frac{(\alpha+\beta)!}{(\alpha+\beta-\gamma)!} \, 
(\xi-\p)^{\beta-\gamma} (\xi^\alpha \ast z^{\alpha+\beta-\gamma}) \\
&=\sum_{\gamma\in \bN^n} \binom \beta \gamma 
\frac{(\alpha+\beta)!}{(\alpha+\beta-\gamma)!} \, 
\xi^{\beta-\gamma} \ast (\xi^\alpha \ast z^{\alpha+\beta-\gamma}) \\
&= \sum_{\gamma\in \bN^n}^\beta \binom \beta \gamma 
\frac{(\alpha+\beta)!}{(\alpha+\beta-\gamma)!} \, 
(\xi^{\alpha+\beta-\gamma} \ast z^{\alpha+\beta-\gamma}).
\end{align*} }

By switching $\xi\leftrightarrow z$ in the argument above and 
using the commutativity of the product $\ast$,  
it is easy to see that we also have 
\begin{align*}
z^\beta(\xi^{\alpha+\beta}\ast z^\alpha) 
= \sum_{\gamma\in \bN^n}^\beta \binom \beta \gamma 
\frac{(\alpha+\beta)!}{(\alpha+\beta-\gamma)!} \, 
(\xi^{\alpha+\beta-\gamma} \ast z^{\alpha+\beta-\gamma}).
\end{align*}
Hence Eq.\,(\ref{even-e1}) follows. 
\epfv

\subsection{Some Applications of Theorem \ref{ast2LP}}
\label{S5.3}
First, let us derive some identities for 
the exponential series 
$\exp_*(\cdot)=e_\ast^{\{\cdot\}}$ of 
the algebra $\cB[\xiz]$, i.e.\@ 
the usual exponential series 
but with the product replaced by 
$\ast$.

\begin{propo}\label{Exp-identitites}
Let $u=(u_1, u_2, ..., u_n)$ be $n$ free commutative variables. 
Set $\xi\ast z\!:=(\xi_1\ast z_1, \, 
\xi_2\ast z_2, ..., \, \xi_n\ast z_n)$ and 
$(\xi\ast z)u\!:=\sum_{i=1}^n (\xi_i\ast z_i)u_i$.
Then, for any ${\bf k}=(k_1, k_2, ..., k_n)\in \bN^n$, we have
\begin{align} 
\xi^{-\bf k}\left(\xi^{\bf k} \ast e_*^{-(\xi \ast z)u}\right) 
&=\prod_{i=1}^n \,\,
\frac{ \exp(-\, 
\frac{(\xi_i z_i) \, u_i}{1-u_i})}{ (1-u_i)^{k_i+1} }. 
\label{Exp-identitites-e1} \\ 
z^{-\bf k}\left(z^{\bf k} \ast e_*^{-(\xi \ast z)u}\right)
&=\prod_{i=1}^n \,\,
\frac{ \exp(-\, 
\frac{(\xi_i z_i) \, u_i}{1-u_i})}{ (1-u_i)^{k_i+1} }. 
\label{Exp-identitites-e2}
\end{align}
In particular, when ${\bf k}=0$, we have the following expression of the exponential $\exp_*(-(\xi \ast z)u)$:
\begin{align}\label{Exp-identitites-e3}
\exp_*(-(\xi \ast z)u)
&=\prod_{i=1}^n \,\,
\frac{ \exp(-\, 
\frac{(\xi_i z_i) \, u_i}{1-u_i})}{ (1-u_i)}. 
\end{align}
\end{propo}

\pf We give a proof for Eq.\,(\ref{Exp-identitites-e1}). 
The proof of Eq.\,(\ref{Exp-identitites-e2}) is similar.

First, by the commutativity and 
associativity of the product 
$\ast$ and also by Proposition \ref{freeness}, 
$(b)$, it is easy to see that, 
for any $\alpha, \beta\in \bN^n$, 
we have 
\begin{align}\label{Exp-identitites-pe1}
(\xi^\alpha \ast z^\alpha)\ast 
(\xi^\beta \ast z^\beta)
&=\xi^{\alpha+\beta} \ast z^{\alpha+\beta}.\\
(\xi \ast z)^{\ast\alpha}&= (\xi^\alpha \ast z^\alpha),
\end{align}
where $(\xi \ast z)^{\ast\alpha}$ denotes the 
``$\alpha^{\rm th}$" power of $(\xi \ast z)$ 
with respect to the new product $\ast$.

By the last two equations above and 
Eq.\,(\ref{ast2LP-e1}), we have
\begin{align}\label{Exp-identitites-pe2}
\xi^{-{\bf k}}\left(\xi^{\bf k}
 \ast \exp_*^{-(\xi \ast z)u}\right) 
&=\sum_{\alpha\in \bN^n} 
\frac{ (-1)^{|\alpha|} }{\alpha!} \, 
\xi^{-\bf k}(\xi^{\alpha+{\bf k}} \ast z^\alpha) 
u^\alpha \\
&=\sum_{\alpha\in \bN^n} 
L_\alpha^{[\bf k]}(\xi z) u^\alpha. \nno
\end{align}

On the other hand, by Eqs.\,(\ref{Generatn-Fctn}) and
(\ref{Mult-LPs}), 
we see that the  generating function of 
the multi-variable generalized Laguerre polynomials 
$L_\alpha^{[\bf k]}(z)$ 
$(\alpha\in \bN^n)$ is given by 
\begin{align}\label{Exp-identitites-pe3}
\prod_{i=1}^n \,\,
\frac{ \exp(-\, 
\frac{z_i \, u_i}{1-u_i})}{ (1-u_i)^{k_i+1} }
=\sum_{\alpha\in \bN^n} 
L_\alpha^{[\bf k]}(z) u^\alpha. 
\end{align}

Replacing $z$ by $\xi z$ in the equation above, 
we get 
\begin{align}\label{Exp-identitites-pe4}
\prod_{i=1}^n \,\,
\frac{ \exp(-\, 
\frac{(\xi_i z_i) \, u_i}{1-u_i})}{ (1-u_i)^{k_i+1} }
=\sum_{\alpha\in \bN^n} 
L_\alpha^{[\bf k]}(\xi z) u^\alpha. 
\end{align}

Combining Eqs.\,(\ref{Exp-identitites-pe2}) 
and (\ref{Exp-identitites-pe4}), we get 
Eq.\,(\ref{Exp-identitites-e1}). 
\epfv

Next we use the connection given in Theorem \ref{ast2LP} to 
derive more properties on the monomials in 
$\xi$ and $z$ with respect to the product $\ast$ 
from certain results on the generalized Laguerre polynomials. 

For convenience, for any $\alpha \in \bN^n$, we set
\begin{align}\label{Def-L-xi}
L_\alpha (z; \xi)\!:=\xi^\alpha \ast z^\alpha.
\end{align}

Note that, by Eqs.\,(\ref{Explicit-Formula}) 
and (\ref{ast2LP-e3}), the polynomials 
$L_\alpha (z; \xi)$ 
$(\alpha\in \bN^n)$ are polynomials 
with coefficients in $\bQ$. In particular, 
for any fixed  $\xi \in (\bR_{>0})^n$, 
by Eqs.\,(\ref{Explicit-Formula}) 
and (\ref{ast2LP-e3}),  
it is easy to see that 
the polynomials $L_\alpha (z; \xi)$ 
$(\alpha\in \bN^n)$ are polynomials 
in $z$ with real coefficients 
and form a linear basis 
of $\cz$. 

The next proposition says that
this basis is also orthogonal  
with respect to the following 
{\it weight} function:
\begin{align}\label{xi-wt}
w_\xi(z)\!:=e^{-\langle \xi, z \rangle} 
\prod_{i=1}^n \xi_i,
\end{align}

\begin{propo}\label{xi-ortho}
For any $\alpha, \beta \in \bN^n$, we have 
\begin{align}\label{xi-ortho-e1}
\int_{(\bR_{>0})^n} 
L_\alpha(z; \xi) L_\beta (z; \xi) 
\, w_\xi(z)\, dz 
=(\alpha!)^2\delta_{\alpha, \beta}.
\end{align}
\end{propo}

\pf Note that, under the change of variables 
$z_i \to \xi_i z_i$ $(\lin)$, 
by Eqs.\,(\ref{ast2LP-e3}) and (\ref{Def-L-xi})   
the Laguerre polynomials $L_\alpha(z)$ 
will be changed to  
\begin{align} 
L_\alpha(z) \to 
L_\alpha (\xi z)=\frac { (-1)^{|\alpha|} }{\alpha!} 
L_\alpha (z; \xi).
\end{align}

By Eq.\,(\ref{xi-wt}) and also Eq.\,(\ref{Lgrwt}) with 
${\bf k}=0$, the weight function $w(z)$ of 
the Laguerre polynomials is changed to 
\begin{align} 
w(z)\to w_\xi(z) \prod_{i=1}^n \xi_i^{-1}.
\end{align}

Now, apply the same changing of the variables 
to the integral in Eq.\,(\ref{ortho-e1}) 
with ${\bf k}=0$, by the last two equations 
above, we get 
\begin{align}
\delta_{\alpha, \beta}
&= \frac { (-1)^{|\alpha+\beta|} }{\alpha!\beta !}
\int_{(\bR_{>0})^n} 
L_\alpha (z; \xi ) L_\beta (z; \xi) 
\, w_\xi(z) \, dz
\end{align}
Hence Eq.\,(\ref{xi-ortho-e1}) 
follows.
\epfv

Denote by $\cA_\bQ[\xiz]$ the polynomial 
algebra in $\xi$ and $z$ over $\bQ$. 
Next we assume $n=1$ and consider 
the irreducibility of the polynomial 
$L_\alpha (z; \xi)$ $(\alpha\in \bN^n)$ 
as elements of $\cA_\bQ[\xiz]$. 
But, first, we need to prove 
the following lemma.

\begin{lemma}\label{Irreduce-lemma}
Let $\xi$ and $z$ be two commutative free 
variables and $K$ any field.
Then, for any $f(z)\in K[z]$ with 
$\deg f\ge 2$,  $f(z)$ is  
irreducible over $K$ iff 
$f(\xi z)\in K[\xi, z]$ $($as a polynomial 
in two variables$)$
is irreducible over $K$. 
\end{lemma}

\pf The $(\Leftarrow)$ part of the lemma is trivial. 
We use the contradiction method to show the 
$(\Rightarrow)$ part of the lemma.

Assume that $f(\xi z)$ is reducible in $K[\xiz]$. Write    
\begin{align}\label{Irreduce-lemma-pe1}
f(\xi z)=g(\xiz) h(\xiz)
\end{align}
for some $g(\xiz), h(\xiz)\in K[\xiz]$
with $\deg g, \deg h\ge 1$. 

Setting $\xi=1$ in the equation above, we also have
\begin{align}\label{Irreduce-lemma-pe1-2}
f(z)=g(1, z) h(1, z).
\end{align}

Let $\bar K$ be the algebraic closure of $K$. Write 
$f(z)=b\prod_{i=1}^d (z-a_i)$ for some $b\in K\backslash \{0\}$ 
and $a_i\in \bar K$ $(1\le i\le d)$. Then we have 
\begin{align}\label{Irreduce-lemma-pe2}
f(\xi z)=b\prod_{i=1}^d (\xi z-a_i). 
\end{align}

Since $f(z)$ is irreducible over 
$K$ and $\deg f\ge 2$ by the assumption, we have  
$a_i\ne 0$ $(1\le i\le d)$. 
Hence, for each $i$, 
$\xi z-a_i$ is irreducible 
in $\bar K[\xiz]$.
Then by Eqs.\,(\ref{Irreduce-lemma-pe1}) 
and (\ref{Irreduce-lemma-pe2}), we have
\begin{align}\label{Irreduce-lemma-pe3}
g(\xi, z)=c\prod_{k=1}^m (\xi z-a_{i_k}) 
\end{align}
for some $c\in \bar K \backslash \{0\}$, $1\le m<d$ 
and $1\le i_1<i_2<\cdots <i_m \le d$. 

However, the equation above implies 
$g(1, z)=c\prod_{k=1}^m (z-a_{i_k})$. 
Since $g(\xi, z)\in K[\xiz]$, we also have  
$g(1, z)\in K[z]$. Then by Eq.\,(\ref{Irreduce-lemma-pe1-2}), 
$g(1, z)$ is a divisor of $f(z)$ in $K[z]$ with 
$1\le \deg g(1, z)=m<d=\deg f(z)$, which 
contradicts to the assumption that $f(z)$ 
is irreducible in $K[z]$.
\epfv 

\begin{theo}\label{Irreduce}
Let $\xi$ and $z$ be two 
commutative free 
variables. For any $m\ge 2$, 
$L_m(z; \xi)=\xi^m \ast z^m$ is 
irreducible in $\cA_\bQ[\xiz]$. 
\end{theo}

\pf By a theorem proved by I. Schur 
\cite{S1}, we know that, for any $m\ge 1$, 
the Laguerre polynomials $L_m(z)$ in one variable 
is irreducible over $\bQ$. Hence, 
by Eq.\,(\ref{ast2LP-e3}) and 
Lemma \ref{Irreduce-lemma}, 
the theorem holds.
\epfv

Note that I. Schur also proved in \cite{S2} 
that the generalized Laguerre 
polynomials $L_m^{[1]}(z)$ $(m\ge 0)$
in one variable  
are also irreducible over $\bQ$. 
Furthermore, M. Filaseta and 
T.-Y. Lam proved in \cite{FL}
that, for any non-negative 
$k\in \bQ$,  all but finitely many of the 
 generalized Laguerre polynomials 
$L_m^{[k]}(z)$ $(m\ge 0)$
in one variable  
are irreducible over $\bQ$. 
Hence, by a similar argument 
as for Theorem \ref{Irreduce}, 
we also have the following 
theorem.

\begin{theo}\label{Irreduce-2}
Let $\xi$ and $z$ be two 
commutative free variables. 
Then, for any $k \in \bN$, 
all but only finitely many of 
the polynomials  
$z^{-k}(\xi^m\ast z^{m+k})$ 
and $\xi^{-k}(\xi^{m+k}\ast z^m)$ 
$(m\in \bN)$ are irreducible 
over $\bQ$.
\end{theo}

Next, we re-prove some important properties of the generalized Laguerre polynomials by using their expressions given in Theorem \ref{ast2LP}. For simplicity, we here only consider 
the one-variable case. Similar results for the 
multi-variable generalized Laguerre polynomials 
can be simply derived from the one-variable case 
via Eq.\,(\ref{Mult-LPs}).

First, let us look at the following recurrent formulas 
of the Laguerre polynomials in one variable.

\begin{propo}\label{Recur-Formula}
For any $m\ge 1$, we have 
\begin{align}
(m+1)L_{m+1}(z)&=(2m+1-z)L_m(z)-m L_{m-1}(z),
\label{RecurRelation}\\
z L'_m(z)&=m (L_m(z)-L_{m-1}(z))\label{D-RecurRelation}
\end{align}
\end{propo}

\pf Note first that, for any $m\ge 1$, 
by Eqs.\,(\ref{B-Propo-e4}) and (\ref{B-Propo-e5}), 
we have
\begin{align*}
\xi \ast z^m &= (\xi-\p)z^m =\xi z^m-mz^{m-1}, \\
z\ast \xi^m &= (z-\dlt)\xi^m=z\xi^m-m\xi^{m-1}. 
\end{align*}
Hence, we also have 
\begin{align*}
\xi \ast z&=\xi z-1, \\
\xi z^m&=\xi \ast z^m  + m z^{m-1}, \\
z\xi^m&=z \ast \xi^m +m\xi^{m-1}. 
\end{align*}

By the last three equations above 
and also Eq.\,(\ref{Symmetry-e1}), 
we have  
\begin{align*}
(\xi z-1)(\xi^m\ast z^m)&=
(\xi\ast z) (\xi^m\ast z^m)=(z\xi^m)\ast(\xi z^m)\\
&=(z \ast \xi^m+m\xi^{m-1})\ast (\xi\ast z^m+m z^{m-1})\\
&=\xi^{m+1}\ast z^{m+1}+2m \, \xi^m\ast z^m+m^2\xi^{m-1}\ast z^{m-1}.
\end{align*}
Multiply $(-1)^m/m!$ to the equation above and then apply 
Eq.\,(\ref{ast2LP-e3}), we get 
\begin{align*}
(\xi z-1)L_m(\xi z)=-(m+1)L_{m+1}(\xi z)+2m L_m(\xi z)
-m L_{m-1}(\xi z).
\end{align*}
Replace $\xi z$ by $z$ in the equation above, we get
\begin{align*}
(z-1)L_m(z)=-(m+1)L_{m+1}(z)+2m L_m(z)
-m L_{m-1}(z),
\end{align*}
Hence Eq.\,(\ref{RecurRelation}) follows. 

To show Eq.\,(\ref{D-RecurRelation}), 
by Eqs.\,(\ref{B-Propo-e5}) and (\ref{ast-Leibniz-e2}), 
we have, 
\begin{align*}
\xi^m\ast z^m &=z\ast(\xi^m\ast z^{m-1})
=(z-\dlt) (\xi^m\ast z^{m-1}) \\
&=z(\xi^m\ast z^{m-1})-m (\xi^{m-1}\ast z^{m-1})\\
&=\frac 1m z\p(\xi^m\ast z^m)-m (\xi^{m-1}\ast z^{m-1}).
\end{align*}
Multiply $(-1)^m/m!$ to the equation above and then apply 
Eq.\,(\ref{ast2LP-e3}), we get 
\begin{align*}
L_m(\xi z) = \frac 1m z\p (L_m(\xi z))+L_{m-1}(\xi z)=
\frac 1m \xi z L'_m(\xi z)+L_{m-1}(\xi z).
\end{align*}
Replace $\xi z$ by $z$ in the equation above, 
we get
\begin{align*}
L_m(z) =
\frac 1m  z L'_m(z)+L_{m-1}(z).
\end{align*}
Hence Eq.\,(\ref{D-RecurRelation}) follows. 
\epfv

Next, we give a new proof for the following 
important property of the generalized 
Laguerre polynomials in one variable.

\begin{theo}\label{Aso-Lag-Equation}
For any $k, m\in \bN$, $L_m^{[k]}(z)$ solves 
the following so-called 
{\it associated Laguerre differential equation}:
\begin{align}\label{Aso-Lag-Equation-e1}
zf''(z)+(k+1-z)f'(z)+m f(z)=0.
\end{align}
\end{theo}

\pf First, by Eq.\,(\ref{Explicit-Formula}), 
we have $L_0^{[k]}(z)=1$. It is easy to see  
that the theorem holds for this case. 

Assume $m\ge 1$. Then, by Eq.\,(\ref{B-Propo-e5}), 
we have 
\begin{align*}
\xi(\xi^{m+k}\ast z^m) 
&=\xi(z\ast(\xi^{m+k}\ast z^{m-1})) \\
&=\xi(z-\dlt) (\xi^{m+k}\ast z^{m-1}) \\
&=\xi z (\xi^{m+k}\ast z^{m-1}) 
-\xi\dlt(\xi^{m+k}\ast z^{m-1}).
\end{align*}
Add $z\p(\xi^{m+k}\ast z^{m-1})$ to the equation above 
and apply Eq.\,(\ref{monomial-eta-e1}), we have 
\begin{align*}
&\xi(\xi^{m+k}\ast z^m)+ z\p(\xi^{m+k}\ast z^{m-1}) \\
&=\xi z (\xi^{m+k}\ast z^{m-1})
-(\xi\dlt-z\p)(\xi^{m+k}\ast z^{m-1})\\
&=(\xi z-k-1)(\xi^{m+k}\ast z^{m-1}).
\end{align*}

By Eq.\,(\ref{ast-Leibniz-e1}), we may re-write the equation above as 
\begin{align*}
\xi (\xi^{m+k}\ast z^m)+ \frac 1m z\p^2(\xi^{m+k}\ast z^{m}) 
=\frac 1m (\xi z-k-1) \p (\xi^{m+k}\ast z^m).
\end{align*}
Multiply $\frac{(-1)^m\xi^{-k-1}}{(m-1)!}$ to  
both sides of the equation above and 
then apply Eq.\,(\ref{ast2LP-e1}), we have
\begin{align*}
m L_m^{[k]}(\xi z)+ z \xi^{-1} \p^2 (L_m^{[k]}(\xi z))  
= (\xi z-k-1) \xi^{-1} \p (L_m^{[k]}(\xi z)).
\end{align*}
By the Chain rule, the equation above is same as 
\begin{align*}
m L_m^{[k]}(\xi z)+ z \xi (\p^2 L_m^{[k]}) (\xi z) 
= (\xi z-k-1) (\p L_m^{[k]})(\xi z).
\end{align*}
Replace $\xi z$ by $z$, or $z$ by $\xi^{-1}z$ 
in the equation above, we get
\begin{align*}
m L_m^{[k]}(z)+ z \p^2 L_m^{[k]} (z) 
= (z-k-1) \p L_m^{[k]}(z).
\end{align*}
Hence we have proved the theorem.
\epfv

Finally, let us point out the following conjecture 
on the generalized Laguerre polynomials, 
which is a special case of  
Conjecture $3.5$ in \cite{GIC} for 
all the classical orthogonal 
polynomials.

\begin{conj}\label{Conj-Lgr}
For any ${\bf k}\in \bN^n$, the subspace $\cM$ of the polynomial algebra $\cA[z]$ spanned by the generalized Laguerre polynomials 
$L^{[\bf k]}_\alpha(z)$ $(0\ne \alpha\in \bN^n)$  is a Mathieu subspace of $\cA[z]$.
\end{conj}

Despite the vast amount of known results  
on the  generalized Laguerre polynomials in the 
literature, the conjecture above is even still 
open for the classical Laguerre polynomials, 
(i.e.\@ the case with ${\bf k}=0$)  
in one variable.

{\small \sc Department of Mathematics, Illinois State University,
Normal, IL 61790-4520.}

{\em E-mail}: wzhao@ilstu.edu.

\end{document}